\documentclass[12pt]{amsart}
\usepackage{amsfonts}
\usepackage{amsmath}
\usepackage{amsxtra}
\usepackage{amssymb,latexsym}
\usepackage[mathcal]{eucal}
\usepackage{amscd}

\makeindex

\newtheorem{theo}{{\bfseries Theorem}}[section]
\newtheorem{prop}[theo]{{\bfseries Proposition}}
\newtheorem{lem}[theo]{{\bfseries Lemma}}
\newtheorem{cor}[theo]{{\bfseries Corollary}}
\newtheorem{df}[theo]{{\bfseries Definition}}

\newtheorem{ex}[theo]{{\bfseries Example}}

\def \N {\mathbb N}

\def \Z {\mathbb Z}
\def \R {\mathbb R}
\def \A {\mathcal A}

\def \CC {\mathcal C}

\def \O {\mathcal O}

\def \ep {\epsilon}
\def \g {\gamma}
\def \l {\lambda}
\def \o {\omega}
\def \d {\delta}
\def \s {\sigma}

\def \th {\theta}

\def \tt {{\mathbf t}}
\def \ss {{\mathbf s}}

\def \zz {{\mathbf z}}

\usepackage{amssymb,latexsym}
\usepackage[mathcal]{eucal}
\usepackage{amscd}
\usepackage{amsmath}
\usepackage{graphicx}
\usepackage{amsfonts}
\usepackage{amscd}

\numberwithin{equation}{section}

\begin{document}

\title{\bfseries Approximation Dynamics}
\vspace{1cm}
\author{Ethan Akin}
\address{Mathematics Department \\
    The City College \\
    137 Street and Convent Avenue \\
    New York City, NY 10031, USA}
\email{ethanakin@earthlink.net}

\date{May, 2019}

\begin{abstract} We describe the approximation of a continuous dynamical system on a p. l. manifold or Cantor set by
a tractable system. A system is tractable when it has a finite number of chain components and, with respect to a given
full background measure, almost every point is generic for one of a finite number of ergodic invariant
measures. The approximations use non-degenerate simplicial dynamical systems for p. l. manifolds
and shift-like dynamical systems for Cantor Sets.

\end{abstract}



\vspace{.5cm} \maketitle

\tableofcontents

\textbf{For Roy Adler, who taught us all so much.}
\vspace{.5cm}

\section{Tractable Dynamical Systems and Approximation}
\vspace{.5cm}

Our spaces will all be compact metric spaces equipped with a metric $d$. For a subset $A$
we use $\overline{A}$ and $A^{\circ}$ for the closure and interior of $A$, respectively.
With $\ep \geq 0$ we let  $\bar V_{\ep} = \{ (x,y) : d(x,y) \leq \ep \}$.\index{$\bar V_{\ep}$}
We let $\Z_+$ \index{$\Z_+$} the set of non-negative integers.
%Our maps will be assumed to be continuous unless otherwise mentioned.
A  \emph{dynamical system} is a pair $(X,f)$ with $f$ a continuous map on a such a space $X$.
The dynamics is given by iterating the map. A subset $A$ of $X$ is called $f$
\emph{invariant} when $f(A) = A$.

Given $\ep > 0$ an $\ep$ \emph{chain}\index{chain} is a sequence $\{ x_0, \dots, x_n \}$
with $n \geq 1$, such that $d(f(x_i),x_{i+1}) \leq \ep $ for
 $0 \leq i < n $. An infinite $\ep$ chain for a map $f$ is often called an $\ep$
 \emph{pseudo-orbit} of $f$. The chain relation $\CC f \subset X \times X$ \index{$\CC f$}
 consists
 of those pairs $(x,y) \in X \times X$ such that for every $\ep > 0$ there exists an
 $\ep$ chain with $x_0 = x$ and $x_n = y$.

 In observing a system, for example in a computer simulation, there is usually some
  positive $\ep$ level of observational
 or computational error. Thus, in attempting to compute an orbit sequence, we are, in fact, producing
 an $\ep$ pseudo-orbit.

A point $x \in X$ is \emph{chain recurrent}\index{chain recurrence} if $(x,x) \in \CC f$ and two such points
$x, y$  are \emph{chain equivalent}\index{chain equivalence}  if
if $(x,y), (y,x) \in \CC f$. The chain equivalence classes of chain recurrent points
are closed, $f$ invariant subsets of $X$ called
\emph{chain components}\index{chain component} or $f$ \emph{basic sets}\index{basic set} (we will use the latter term).

Chain recurrence is the broadest notion of recurrence. If $x = f(x)$ then $x$ is an
\emph{equilibrium point} or \emph{fixed point} for $f$. If $x = f^n(x)$ for some positive integer $n$ then $x$ is a
\emph{periodic point}. Let $\o f(x)$ \index{$\o f$}
denote the set of limit points of the orbit sequence
$\O f(x) = \{ f(x), f^2(x), \dots\}$. If $x \in \o f(x)$, then $x$ is a \emph{recurrent point}. For any $x \in X$, the set
$\o f(x)$ is a nonempty, closed, $f$ invariant subset of $X$ which is contained in a single basic set.

The simplest sort of dynamical system has a finite chain recurrent set.
In that case, each basic set is a periodic orbit and
$\o f(x)$ is a periodic orbit for every $x \in X$. From any initial
point the orbit exhibits a transient motion towards
its limiting periodic orbit. Eventually, we observe what appears to be the
periodic motion of the limit orbit.

We will focus upon the more general case when the chain recurrent set decomposes into only finitely many basic sets.
For a basic set $B$ we will define the \emph{basin}\index{basin} of $B$ to be the set $Bas(B) =
(\o f)^{-1}(B)  = \{ x : \o f(x) \subset B \}$\index{$Bas(B)$}.
We will call $B$ a
\emph{visible basic set}\index{basic set!visible}
 if $Bas(B)$ has a nonempty interior. When there are
only finitely many basic sets then the union of the basin interiors  is a dense
open subset of $X$ (see Theorem \ref{tractheo00a} below).

In that case, the orbit every point approaches a basic set and most points approach
a
% $\mu_0$
visible basic set. However, the
motion within a basic set
may be quite complicated, i.e. chaotic. So to describe such motion we resort to
statistics using an invariant measure.

By a \emph{measure} $\mu$ on $X$ we mean a Borel probability measure. The
\emph{support}\index{measure!support} $|\mu|$ of the measure consists
of those points whose every neighborhood has positive measure.  The measure is
 called  \emph{full }\index{measure!full} when $|\mu| = X$.
The measure is $f$ \emph{invariant}\index{measure!invariant} when $\mu(f^{-1}(A)) = \mu(A)$ for every
Borel set $A$, or, equivalently, when
$\int  u(f(x)) \  \mu(dx) \ = \  \int  u(x) \ \mu(dx)$ for every continuous
real-valued function $u$ on $X$. An invariant measure is
\emph{ergodic}\index{measure!ergodic} if for any Borel set $A$, $f^{-1}(A) = A$ implies $\mu(A) = 0$ or $1$.

A point $x$ is called \emph{generic}\index{generic point} for an invariant measure $\mu$ if
\begin{equation}\label{intro01}
\lim_{n \to \infty} \ \frac{1}{n} \ \Sigma_{i=0}^{n-1} \ u(f^i(x)) \ = \ \int u \  d \mu
\end{equation}
for every continuous
real-valued function $u$. Since the measure is determined by the sequence of
points of the orbit, it is clear that a point is generic for
at most one measure.

We let $Gen(\mu)$ denote the set of points generic for $\mu$. While it need not be true that a generic point lies in $|\mu|$,
the Birkhoff Pointwise Ergodic Theorem says that if $\mu$ is an ergodic measure,
then $\mu(|\mu| \cap Gen(\mu)) = 1$.
If $x \in |\mu| \cap Gen(\mu)$, then $\o f(x) = |\mu|$.  In any case, if
$x \in Gen(\mu)$, then $|\mu| \subset \o f (x)$. So if $Gen(\mu)$ is
 nonempty, e.g. if $\mu$ is ergodic, then $|\mu|$ is contained in a basic set which we then denote $B(\mu)$.

The tractable systems  occur on spaces equipped with some given background measure
$\mu_0$. We will call $(X,f,\mu_0)$ a \emph{measured dynamical system}\index{dynamical system!measured} when
$(X,f)$ is a dynamical system, i.e. $f$ is a continuous map on a compact metric
space $X$, and $\mu_0$ is a full measure on $X$.
The measure $\mu_0$ is usually not invariant, but we assume that it is complete,
i.e. every subset of a set of measure zero is measurable and so has
measure zero.

For a measured dynamical system, we call a basic set $B$ $\mu_0$-visible\index{basic set!$\mu_0$-visible} when
$\mu_0(Bas(B)) > 0$. Since $\mu_0$ is full, a visible basic set
is $\mu_0$-visible.

\begin{df}\label{introdef01} Let $(X,f,\mu_0)$ be a measured dynamical system.
 We will call  the system \emph{tractable}\index{dynamical system!tractable} when

\begin{itemize}

\item[TRAC 1] There are a finite number of ergodic invariant
probability measures $\mu_1, \dots, \mu_k $ such that $\mu_0(Gen(\mu_i)) > 0$ for $i = 1, \dots k$ and
 $\mu_0(\bigcup_{i=1}^{k} \ Gen(\mu_i)) = 1$. We call these the \emph{special ergodic measures} for $f$.

 \item[TRAC 2] There are only finitely many basic sets.

\item[TRAC 3] There are a finite number of pairwise disjoint open subsets of $X$
 $U_1, \dots, U_k $ such that  $\mu_0(Gen(\mu_i) \setminus U_i) = 0$ for $i = 1, \dots k$.
Hence, the union $U = \bigcup_{i=1}^k U_k$ has $\mu_0$ measure one and so is dense in $X$.
%
%\item[TRAC 3] For $i = 1, \dots k$, if $B_i$ is the basic set containing $|\mu_i|$ then $\{ B_1, \dots , B_k \}$
%are the visible basic sets.
%
%\item[TRAC 4] It may happen that a visible basic set contains more than one $|\mu_i|$,
%but if $i \not= j$ then the intersection $|\mu_i| \cap |\mu_j|$ has $\mu_0$ measure zero and so is nowhere dense.
% in $|\mu_i|$.
\end{itemize}
\end{df}
\vspace{.5cm}

%We can in TRAC 2 obtain the condition $\mu_0(Gen(\mu_i)) > 0$ for $i = 1, \dots k$ from the assumption  $\mu_0(\bigcup_{i=1}^{k} \ Gen(\mu_i)) = 1$,
%by discarding those measures with $\mu_0(Gen(\mu_i)) = 0$.
%From TRAC 2 and TRAC 3, the union $U = \bigcup_{i=1}^k U_k$ has measure one and so is dense in $X$.

Thus, for a tractable system, an initial point $x \in X$ is, with $\mu_0$ probability
one, generic for a special ergodic measure $\mu_i$. In particular,
a point $x$ of the region $U_i$ is, almost surely, generic for $\mu_i$.
The orbit sequence then approaches the set $\o f(x) \subset B(\mu_i)$. Furthermore, if we
regard a continuous real-valued function $u$ on $X$ as a
measurement, then the average of the successive measurements along the orbit,
$\frac{1}{n}  \Sigma_{i=0}^{n-1}  u(f^i(x))$, approaches the $\mu_i$ expected
value of $u$ as $n \to \infty$.

Tractable systems do exist, e.g. Morse-Smale diffeomorphisms and the Axiom A
diffeomorphisms described by Smale in \cite{S67}.
On the other hand, the generic homeomorphisms on manifolds of
positive dimension and on Cantor sets have infinitely many attractors and
uncountably many basic sets, see \cite{AHK}. So such a system is certainly not tractable.

Rather than look for conditions which will assure that a
system is tractable, our approach here will be to approximate an arbitrary system by a
tractable ones in a natural way.

Our approximations will be built using finite data. The analogy is with curve sketching. To sketch the graph
of a real-valued function on a bounded interval we plot a finite number of points and
then connect the dots with line segments.
The resulting graph is piecewise linear and is completely determined by
the end-points of each segment, but we use the p.l. function instead of the dots
because then our approximation is
the same sort of object as that which we wish to study, i.e. a real-valued function on the interval.

The advantage of these approximations is that the dynamics is directly
computable from the finite initial data.
The disadvantages are two-fold. First, we are operating in the topological
category so that
our approximating functions $g$ are merely $C^0$ close to the original
function $f$. The $g$ orbit of a point
is a pseudo-orbit of the original function $f$. Second, our functions $g$
are never invertible and so
if $f$ is a homeomorphism and we want to describe the backward as well as
forward dynamics, then we require
separate approximations for $f$ and $f^{-1}$. On the other hand, computer
simulations usually produce only
pseudo-orbits and so we expect that our approximations should pick up the
sorts of things observed by a simulation.

In the next two sections we review the dynamics of  relations. In addition to
the dynamics of a
continuous map on a general compact metric space, we are interested in the case
when $G \subset K \times K$ with $K$ a
finite set. The two situations  are related in that associated with
such a relation $G$ is the \emph{sample path space}, a closed subset $K_G$ of
$K^{\Z_+}$ on which the shift map restricts to a subshift of finite type. A
\emph{stochastic cover}\index{stochastic cover} $\Gamma$ of $G$
is a stochastic matrix with $\Gamma_{ji} > 0$ if and only if $(i,j) \in G$. Such a stochastic
cover induces a Markov chain with associated
Markov measures on $K_G$.

We will see in Section 3 that from the stochastic cover we obtain a background measure so that
the shift on $K_G$ becomes the paradigmatic example of a tractable system. We
also consider  the so-called
\emph{special two-alphabet model} which is applied in the later sections. Such a
model consists of two finite sets $K^*$ and $K$
together with maps $J, \g : K^* \to K$. We then obtain the relations $G = \g \circ J^{-1}$
 and $G^* = J^{-1} \circ \g$ on
$K$ and $K^*$, respectively. The classic example would be a directed graph
with $K^*$ the set of edges and $K$ the set of
vertices. The maps $J$ and $\g$ associate to an edge its initial and terminal
vertex, respectively. When the map $J$ in the
two-alphabet model is surjective, we can choose for each $s \in K$ a distribution
on $J^{-1}(s) \subset K^*$.
From this \emph{distribution
data}\index{distribution data} there naturally arise stochastic covers for the relations $G$ and $G^*$.

In Section 4,  $X$ is a $d$-dimensional p.l. manifold. We
approximate the continuous map $f$ on $X$ by a non-degenerate \emph{simplicial
dynamical system}\index{dynamical system!simplicial}.
Let $L$ is a simplicial complex triangulating $X$ and $L^*$ be a subdivision
which is \emph{proper}\index{proper subdivision}, i.e. no simplex of $L^*$ meets disjoint simplices of $L$. Let
$K$ and $ K^*$ be the sets of $d$-dimensional
simplices in each complex. Let $J :$  $ K^* \to $ $K$ be the
map associating to each $s^* \in$ $K^*$ the unique $s \in$ $K$ which
contains it.  A non-degenerate simplicial dynamical system\index{dynamical system!simplicial!nondegenerate}
is a simplicial map $\g : L^* \to L$ such that the dimension of
$ s = \g(s^*) \in L$ is equal to that of $s^*$ for every simplex  $s^* \in L^*$. So $\g$ restricts to a
map from $K^*$ to $K$. Furthermore, $\g$ can be chosen so that the associated p.l. map $g$ on $X$
is uniformly close to $f$. A simplicial map is determined by its value on the
vertices and so is given by finite data.
 The relation $G^* = J^{-1} \circ \g$ on the finite set  $ K^*$ induces the
 subshift   $(K^*_{G^*},S)$.
There is an almost one-to-one map taking subshift onto $(X,g)$, thus providing a
simple description of the dynamics of $g$.
Furthermore, with
a locally Lebesgue measure $\l_0$ taken as background measure on $X$, the system
$(X,g, \l_0)$ is tractable.

Finally, in Section 5 $X$ is a Cantor set $A^{\Z_+}$ with $A$ a finite alphabet.
For $z \in X$ and $s$ a finite word in the alphabet $A$,
we write $x = sz \in X$ for the concatenation of $s$ with $z$.  Fix two positive integers $n$ and $k$.
Let $K$ and $K^*$ be the sets of words in the alphabet
of length $n$ and $n+k$, respectively. Define $J_n : X \to K$ and $J_{n+k} : X \to K^*$
associating to $x \in X$ the
initial words of length $n$ and $n+k$, respectively. Similarly, $J : K^* \to K$ maps
each $n+k$ word to its initial $n$ word. If $\g : K^* \to K$
is an arbitrary map, the associated \emph{shift-like} map $g$ on $X$ is given by
$g(s^*z) = \g(s^*)z$ where $s^* = J_{n+k}(x)$ for
$x = s^*z$ and so $\g(s^*) = J_{n}(g(x))$.  The relation $G^* = J^{-1} \circ \g$ on the
finite set  $ K^*$ induces the subshift $(K^*_{G^*},S)$. There is a homeomorphism taking the subshift
onto $(X,g)$ and so the two systems are conjugate. Furthermore, with the uniform Bernoulli
 measure $\mu_0$ taken as background measure on $X$, the system
$(X,g,\mu_0)$ is tractable. For any continuous function $f$ on $X$ and any positive integer $n$,
there exists a positive integer $k$
so that $J_n(f(x))$ depends only on $J_{n+k}(x)$ and so $\g = J_n \circ f \circ J_{n+k}^{-1} $
is a map from $K^*$ to $K$.
If $g$ is the shift-like map associated with $\g$ then for all $x \in X$,
$J_n(f(x)) = J_n(g(x))$ and so $g$ approximates
$f$.  In addition, there is a continuous map $Q^f$ on $X$ which maps $f$ to $g$. For
every $x$, $J_{n+k}(Q^f(x)) = J_{n+k}(x)$.
Since $Q^f$ maps $f$ to $g$, the $g$ orbit of $Q^f(x)$ shadows the $f$ orbit of $x$.
\vspace{1cm}

\section{Relation Dynamics And Tractable Systems}
\vspace{.5cm}

One can define many dynamics notions when one merely assumes
that $F$ is a closed relation on a space $X$, i.e. a closed
subset of $X \times X$. We begin by outlining some of the relation properties
defined and described in Chapters 1 and 2 of
\cite{A93}. We describe the general setup as we will need it for  relations on finite
sets as well as for continuous maps on compact metric spaces.

A (closed) \emph{relation}\index{relation}
$F : X \to Y$ is just a (closed) subset of $X \times Y$.  We let
$F(x) = \{ y : (x,y) \in F \}$ and $F(A) = \bigcup_{x \in A} F(x)$. For example,
$\bar V_{\ep}(x)$ is the closed $\ep$ ball centered at $x$.
The set $F(A)$ is the projection to $Y$ of the set $F \cap (A \times Y) \subset X \times Y$.
So if $F$ and $A$ are closed, then $F(A)$ is closed.

For $F : X \to Y$, we define $F^{-1}: Y \to X$ by
$F^{-1} = \{ (y,x) : (x,y) \in F \}$. The \emph{domain}\index{domain} of $F$ is $Dom(F) = F^{-1}(Y) = \{ x : F(x) \not= \emptyset \}$.\index{$Dom(F)$}
 The relation $F$ is called  \emph{surjective}\index{relation!surjective}
$X = F^{-1}(Y)$ and $Y = F(X)$.

The relation $F$ is a map, exactly when $F(x)$ is a singleton for every
$x$. In that case we will use the same symbol $F(x)$ for the set and the point therein.
For example, the identity map $1_X$ on $X$ is the diagonal set $\{ (x,x) : x \in X \}$.
A map is continuous if and only if it is closed as a relation. It is surjective as a map if and only if it is surjective as a relation.

If $F : X \to Y, G : Y \to Z$ are relations then the composition
$G \circ F : X \to Z  = \{ (x,z) :$ for some $y \in Y,  (x,y) \in F$ and $ (y,z) \in G\}$.
Thus, $G \circ F$ is the projection to
$X \times Z$ of  $(F \times Z) \cap (X \times G) \ \subset \ X \times Y \times Z$. So the
composition of closed relations is closed.  Clearly, $(G \circ F)^{-1} = F^{-1} \circ G^{-1}$.

If $F$ is a relation on $X$, i.e. a subset of $X \times X$,
then we define $F^0 = 1_X$, the identity map; $F^1 = F$ and, inductively,
$F^{n+1} = F \circ F^n$.  Since composition is associative it
follows that $F^{n + m} = F^n \circ F^m$ for any non-negative integers $n, m$.
We let $F^{-n} = (F^{-1})^n$.

If $A \subset X$ and $F$ is a relation on $X$, $A$ is called $^+$\emph{invariant}
if $F(A) \subset A$ and \emph{invariant}\index{invariant set} if $F(A) = A$. If $F$ and $A$ are closed
and $A$ is $^+$invariant, then
$\bigcap_{n=0}^{\infty} F^n(A)$ is invariant.

 If $A$ is a closed subset of $X$, then the \emph{restriction} $F|A$ of a closed
 relation $F$ on $X$ is the closed relation
 $F \cap (A \times A)$ on $A$. If $F$ is a map on $X$, then the restriction $F|A$ is
 a map on $A$ when $A$ is $^+$invariant.

A relation $F$ on $X$ is \emph{reflexive} if $1_X \subset F$, \emph{symmetric} if
$F^{-1} = F$ and \emph{transitive} if
$F \circ F \subset F$.

From a closed relation $F$ on $X$ we construct other relations on $X$.

\begin{itemize}
\item The \emph{orbit relation}: $\O F = \bigcup_{n = 1}^{\infty} \ F^n $.\index{orbit relation} \index{$\O F$ }

\item The \emph{limit point relation}: $\o F$, defined by \\
$\o F(x) = \bigcap_{n = 1}^{\infty} \overline{\bigcup_{m \geq n} F^m(x)}.$\index{limit point relation} \index{$\o F$}

%\item The \emph{wandering relation}: $\NN f = \overline{\bigcup_{n = 1}^{\infty} \ f^n} $.

\item The \emph{chain relation}: $\CC F = \bigcap_{\ep > 0} \O (\bar V_{\ep} \circ F)$.\index{chain relation} \index{$\CC F$}
\end{itemize}

 The orbit relation is transitive, but neither $\O F$ nor $\o F$ is usually closed.
The chain relation $\CC F$ is both closed and transitive. Notice that  if $X$ is finite,
and so is discrete, then $\O F =  \CC F$.

\begin{equation}\label{inveq}
\O (F^{-1}) \ = \ (\O F)^{-1}, \qquad  \CC (F^{-1}) \ = \ (\CC F)^{-1}
\end{equation}
and so in these cases we can omit the parentheses (for the latter equation, see \cite{A93} Prop. 1.11(a)).  On the other hand, $\o (F^{-1})$
is usually not equal to $ (\o F)^{-1}$ even when $F$ is a homeomorphic map.

 For a relation $F$ on $X$ we let $|F| =   \{ x \in X : (x,x) \in  F \}$ \index{$|F| $}
 which is a
 closed set when $F$ is a closed relation.
 Extending the language from the case when $F$ is a continuous map,  we call the points
 of $|F|$ the \emph{fixed} points of $F$,\index{fixed point}
 $|\O F|$ the \emph{periodic} points,\index{periodic point}  $|\o F|$ the \emph{recurrent points}\index{recurrent point}
 and $| \CC F |$ the \emph{ chain recurrent}\index{chain recurrent point} points. On the closed set
 $|\CC F|$ the relation $\CC F$  restricts to a reflexive,
 transitive relation and so $\CC F \cap \CC F^{-1}$ is a closed equivalence
 relation on $|\CC F|$.  The equivalence classes are called the $F$
 \emph{basic sets}\index{basic set} or \emph{chain components} of $F$. $\CC F$ induces a partial order on the set of basic subsets.
 We call a closed relation $F$ on $X$  \emph{chain transitive}\index{relation!chain transitive} when $\CC F = X \times X$
 and so entire space $X$ is a single basic set.

 A basic set is called \emph{terminal}\index{basic set!terminal} when it is $\CC F$ invariant. That is,
for a terminal basic set $B$, $x \in B$ implies $\CC F(x) \subset B$. If $Dom(F) = X$, then $B$ is a
terminal basic set if and only if it is minimal in the collection of nonempty $\CC F$ $^+$invariant sets.
In that case, it follows from the usual Zorn's Lemma argument, that any nonempty, closed
$\CC F$ $^+$invariant set contains a terminal basic set.

 Following \cite{A93} Chapter 3, we call  a closed subset $U$ \emph{inward}\index{inward set} for a closed relation $F$ on $X$
when $F(U) \subset U^{\circ}$. An inward set is $\CC F$ $^+$invariant, and  the $\CC F$ invariant set
$$ \bigcap_{n = 0}^{\infty} (\CC F)^n(U) \ = \ \bigcap_{n = 0}^{\infty} F^n(U) \ =  \ A$$
 is called the \emph{attractor}\index{attractor}  to which $U$ is associated. For an attractor $A$,
 the set of inward sets associated to it form a neighborhood
  base for $A$. If $x \in B \cap A$ with $B$ a basic set and $A$ an attractor, then $B \subset \CC F(x) \subset A$.
  Thus, $A$ contains any basic set that it meets.

  If $U$ is inward for $F$, then $U' = X \setminus U^{\circ}$ is
 inward for $F^{-1}$. The attractor $R$ for $F^{-1}$ to which $U'$ is associated is
 called the \emph{repeller dual to $A$}\index{repeller} and $(A,R)$ is called
 an \emph{attractor-repeller pair}\index{attractor-repeller pair} for $F$. A pair of closed sets $(A,R)$ is an
 attractor-repeller pair when it satisfies
 \begin{equation}\label{attrepeq1}
 \begin{split}
 A \ = \ \CC F(A), \qquad R \ = \ \CC F^{-1}(R), \\
 A \cap B \ = \ \emptyset, \qquad |\CC F| \subset A \cup B.
 \end{split}
 \end{equation}
 Furthermore, an attractor-repeller pair is determined by the basic sets each contains.
  \begin{equation}\label{attrepeq2}
 A \ = \ \CC F(A \cap |\CC F|), \qquad  R \ = \ \CC F^{-1}(R \cap |\CC F|).
 \end{equation}
 See \cite{A93} Propositions 3.8, 3.9.

 Now we specialize to the case of a dynamical system $(X,f)$, i.e. the relation $f$ is a continuous map on $X$.

 When $f$ is a continuous map, each basic set and each limit point set $\o f(x)$ is $f$
 invariant. Furthermore, if $(A,R)$ is an attractor-repeller
 pair, then $A$ is $f$ invariant and $R$ is $f$ $^+$invariant, by \cite{A93} Proposition 3.8 and Corollary 4.3.

 As with a general closed relation, a continuous map $f$ on $X$ is chain transitive when
 $\CC f = X \times X$, i.e. the entire space
 $X$ is a basic set for $f$.
 Since a basic set is invariant, it follows that a chain transitive map is surjective.

 A continuous map $f$ on $X$ is \emph{topologically transitive}\index{map!topologically transitive} if for every pair of
 nonempty open subsets $U, V$ of $X$,
 there exists $i \in \Z_+$ such that $U \cap f^{-i}(V) \not= \emptyset $.
We call a point $x \in X$ a \emph{transitive point}
when $\O f(x)$ is dense in $X$, or, equivalently,
when $\o f(x) = X$.
A map is topologically transitive if and only if it admits a transitive point and
in that case the set of transitive points is an invariant, dense $G_{\delta}$ subset of $X$.

If $(X,f)$ is a dynamical system and  $A \subset X$ is closed and $^+$invariant, then $(A,f|A)$ is the
 subsystem obtained by restricting $f$ to $A$. A closed $^+$invariant subset
 $A$ is called  a \emph{chain transitive subset} (or a \emph{topologically transitive subset})
 for a continuous map $f$ when the restriction $f|A$
 is a chain transitive map (resp. a topologically transitive map). For example, for a continuous map $f$
 each limit point set $\o f(x)$ is a chain transitive subset (Proposition 4.14
 of \cite{A93}) and so is contained in a single basic set.

For a  map $f$ and a closed, $^+$invariant set $A$, we will call the Borel set
 \begin{align}\label{attrepeq2a}
\begin{split} Bas(A) \ &=  \{ x : \o f (x) \subset A \} \\ = \
 \bigcap_{k=1}^{\infty} \bigcup_{N=1}^{\infty} &\bigcap_{n=N}^{\infty} f^{-n}(\bar V_{1/k}(A))
 \end{split} \end{align}
 the \emph{basin} of $A$.

\begin{lem}\label{attreplem1} For a dynamical system $(X,f)$ let $A$ be
an attractor with dual repeller $R$.
 \begin{equation}\label{attrepeq3}
 X \setminus R \ = \ \{ x : \o f(x) \cap A \not= \emptyset \} \ = \ Bas(A).
 \end{equation}
 In particular, the basin of an attractor is an open set.
 \end{lem}

 {\bfseries Proof:} Every basic set is entirely contained in $A$ or $R$,
 and any limit point set is contained in a basic set.  It follows that for every $x \in X$, either
 $\o f(x) \subset A$ or $\o f(x) \subset R$. Because $R$ is closed and
 $f$ $^+$invariant, $x \in R$ implies $\o f(x) \subset R$.
On the other hand, if $y \in \o f(x) \cap R$, then
$x \in (\o f)^{-1}(y) \subset \CC f^{-1}(y) \subset R$. Contrapositively,
$x \not\in R$ implies $\o f(x)$ is disjoint from $R$.

$\Box$ \vspace{.5cm}

We will call a basic set  \emph{visible}\index{basic set!visible} when its basin has a non-empty interior.

Recall that a basic set is called \emph{terminal} when it is $\CC f$ invariant. That is,
for a terminal basic set $B$, $x \in B$ implies $\CC f(x) \subset B$. Since $Dom(f) = X$,  $B$ is a
terminal basic set if and only if it is minimal in the collection of nonempty $\CC f$ $^+$invariant sets.

\begin{theo}\label{tractheo00a} Assume that $(X,f)$ is a dynamical system with only finitely many basic sets.
Each terminal basic set is an attractor and so is visible. Furthermore,
the union of the basin interiors of the basic sets is a dense open subset of $X$. \end{theo}

{\bfseries Proof:}  Let $B$ be a terminal basic set. Number the basic sets
$\{ B_1, \dots , B_n \}$ so that $B_n = B$.
We proceed by induction on $n$.

Let $Y = \CC f^{-1}(\bigcup_{i < n} B_i)$. Because $B_n$ is $\CC f$ invariant, $Y$ is
disjoint from $B_n$. It follows from
(\ref{attrepeq1}) that $(B_n,Y)$ is an attractor-repeller pair and so from Lemma \ref{attreplem1} it follows that
$Bas(B_n) = X \setminus Y$.

%Since there are only finitely many basic sets, $B_n \cap |\CC f| = |\CC f| \setminus (\bigcup_{i < n} B_i) \cap |\CC f|$ and
%so is clopen in $|\CC f|$. Since $B_n$ is terminal, $\CC f(B_n) = B_n$ and so $B_n$ is an attractor by Theorem 3.3 of \cite{A93}.
%From Proposition 3.9 of \cite{A93} it follows that $\{ x : \o f (x) \subset B_n \}$ is the basin $G(B_n)$ and its complement,
%$Y$, is the dual repeller.
%This implies that $\CC f^{-1}(Y) = Y$.
Let $\hat G$ be the interior of $Y$.
Since $Y$ is closed, $Bas(B_n) \cup \hat G$ is a dense open
subset of $X$.

If $\hat  G = \emptyset$ then the union of the basin interiors is $Bas(B_n)$ and it is dense in $X$.
In particular, this proves the initial -  $n = 1$ -
step of the induction.

Assume instead that $\hat G$ is non-empty.

By Proposition 2.6 and Corollary 4.3 of \cite{A93}, $Y$ is a closed $^+$invariant subset
and by Theorem 3.5 of \cite{A93} $\CC (f|Y) = (\CC f)|Y$. This implies that the basic sets
of the subsystem $(Y,f|Y)$ are
$B_1, \dots, B_{n-1}$. Let $G_i$ be the $Y$ interior of the $f|Y$ basin of the basic set
$B_i$ for $i < n$.
Each $G_i$ is open in the relative topology of $Y$ so $G_i \cap \hat G$ is open in $X$.  This is clearly contained in
(and, in fact, equals) the basin interior  $Bas(B_i)^{\circ}$ for $(X,f)$. By induction
hypothesis, $\bigcup_{i < n} G_i$ is dense in $Y$.
Since $\hat G$ is open, $\bigcup_{i < n} G_i \cap \hat G$ is dense in $\hat G$. It follows that
$Bas(B_n) \cup (\bigcup_{i < n} G_i \cap \hat G) \subset \bigcup_{i \leq n} Bas(B_i)^{\circ} $ is dense in $X$.
This completes the proof of the inductive step.

$\Box$ \vspace{.5cm}

For every point $x$, the limit point set $\o f(x)$ is contained in some basic set.
Thus, eventually, the motion is
close to some basic set. To analyze the behavior within or close to a basic set, we use invariant measures.

The space of signed measures on $X$ is the dual space to the Banach space of
continuous, real-valued functions
on $X$.  On it we use the weak$^*$ topology and in it the set $P(X)$ of probability
measures is a compact, convex subset.
From now on, we will assume that our measures are probability measures unless otherwise mentioned.
For $x \in X$ we let $\d_x$ denote the point measure concentrated on $x$.
If $\mu \in P(X)$, then the \emph{support} $|\mu| = \{ x : \mu(U) >  0 $ for every
neighborhood $U$ of $ x \}$. That is,
the closed set $|\mu|$ is the complement of the union of the open sets $U$ of $\mu$ measure zero.
Since $X$ has a countable base,
$X \setminus |\mu|$ is the maximum open set of measure zero.  $\mu$ is called \emph{full} if
$|\mu| = X$, i.e. every non-empty open subset has positive measure.

If $f : X_1 \to X_2$ is a continuous map, and $\mu_1 \in P(X_1)$, then $f_*(\mu) \in P(X_2)$ is
defined by $f_*(\mu)(A) = \mu(f^{-1}(A))$ for
all Borel sets $A \subset X_2$ or, equivalently, by $\int u(y) f_*(\mu)(dy) = \int u(f(x)) \mu(dx)$
for every continuous
$u : X_2 \to \R$. For example, $f_*(\d_x) = \d_{f(x)}$.

Thus, the continuous map $f$ induces the continuous map $f_* : P(X_1) \to P(X_2)$.  It is
clear that $f_*$ is affine, i.e.
$f_*(t \mu_1 + (1-t) \mu_2) \ = \ t f_*(\mu_1) + (1-t) f_*(\mu_2).$ for $t \in [0,1]$.

If $f$ is a map on $X$ then
$f_*$ is a map on $P(X)$. We call $\mu$ an \emph{invariant} measure when $f_* \mu = \mu$. That is,
$|f_*| \subset P(X)$ is the compact, convex set
of invariant measures.

For measures $\mu, \nu$ on $X$ we say that $\nu$ is \emph{absolutely continuous} with
respect to $\mu$, written
$\mu \gg \nu$, if $\mu(A) = 0$ implies $\nu(A) = 0$ for all Borel sets $A$. The measures
are \emph{absolutely equivalent},
written $\mu  \approx \nu$, if $\mu \gg \nu$ and $\nu \gg \mu$, i.e. they have the same
sets of measure zero.

A measure $\mu \in |f_*|$ is \emph{ergodic} if $A \subset f^{-1}(A)$ for a Borel set $A$
implies $\mu(A) = 0$ or $1$, or, equivalently,
if $u : X \to \R$ is a bounded measurable function such that $u \circ f = u$ a.e. then $u$
is the constant function
$u(x) = \int u(y) \mu(dy)$ for $\mu$ almost every
$x \in X$.

A point $x$ is called \emph{generic} for $\mu$ if for every continuous $u : X \to \R$,
$\lim_{n \to \infty} \ \frac{1}{n} \Sigma_{i = 0}^{n-1} u(f^i(x)) =
\int  u(x) \mu(dx)$, or, equivalently, if $\mu$ is the limit in $P(X)$ of the
sequence $\{ \frac{1}{n} \Sigma_{i=0}^{n-1} \d_{f^i(x)} \}$. We write $Gen(\mu)$\index{$Gen(\mu)$} for the
(possibly empty) set of points generic for $\mu$.
As with basin of a fixed point, it is clear that
$Gen(\mu)$ is a Borel set.
%If $x \in Gen(\mu)$  and $x_1$ is asymptotic to $x$, i.e. \\ $\lim_{n \to \infty} d(f^n(x),f^{n_1}(x)) = 0$,
%then $x_1 \in Gen(\mu)$.

The \emph{ Birkhoff Pointwise Ergodic Theorem} says that for an ergodic measure $\mu$,
$\mu$ almost every point of $|\mu|$ is generic for $\mu$, i.e. $Gen(\mu) \cap |\mu|$ has $\mu$ measure one.

Notice that we do not assume that $x \in |\mu|$ for $x \in Gen(\mu)$.

\begin{prop}\label{tracprop00aa} For a dynamical system $(X,f)$, let $\mu$ be an invariant measure for $f$.

\begin{itemize}
\item[(a)] The support $|\mu|$ is $f$ invariant.

\item[(b)] A point $x \in Gen(\mu)$ if and only if $f(x) \in Gen(\mu)$, i.e. $Gen(\mu) = f^{-1}(Gen(\mu))$.
In particular, $Gen(\mu)$ is $^+$invariant. It is invariant if $f$ is surjective.
\end{itemize}
\end{prop}

{\bfseries Proof:} (a) Let $x \in X$.

If $x \in |\mu|$ and $U$ is a neighborhood of $f(x)$ then $\mu(U) = \mu(f^{-1}(U)) > 0$ since
$f^{-1}(U)$ is a neighborhood of $x$. Thus, $f(x) \in |\mu|$ and so $|\mu|$ is $f$ $^+$invariant.

Now suppose $f^{-1}(x) \cap |\mu| = \emptyset$. Each point of $f^{-1}(x)$ has a neighborhood
of measure zero. By using the
union of a finite subcover, we obtain an open set $V$ of measure zero which contains $f^{-1}(x)$.
For example, if
$f^{-1}(x) = \emptyset$ then let $V = \emptyset$. There exists an open set $U $ containing
$x$ such that $f^{-1}(U) \subset V$.
Hence, $\mu(U) = \mu(f^{-1}(U)) = 0$. Hence, $x \not\in |\mu|$. Contrapositively, $x \in |\mu|$ implies
$f^{-1}(x) \cap |\mu| \not= \emptyset$. Hence, $|\mu|$ is $f$ invariant.

(b) The sequences $\{ \frac{1}{n} \Sigma_{i=0}^{n-1} \d_{f^i(x)} \}$ and $\{ \frac{1}{n} \Sigma_{i=1}^{n} \d_{f^i(x)} \}$
are asymptotic.

$\Box$ \vspace{.5cm}

\begin{prop}\label{tracprop00ab} For a dynamical system $(X,f)$, let $\mu$ be an invariant
measure for $f$ with $Gen(\mu)$ nonempty..

\begin{itemize}
\item[(a)] There exists a unique basic set denoted $B(\mu)$\index{$B(\mu)$} such that $|\mu| \subset \o f(x) \subset B(\mu)$ for all
$x \in Gen(\mu)$. That is, $Gen(\mu) \subset Bas(B(\mu))$.
If $B$ is a basic set and $Bas(B) \cap Gen(\mu) \not= \emptyset$ or $|\mu| \cap B \not= \emptyset$, then $B = B(\mu)$.

\item[(b)] If $x \in  Gen(\mu) \cap |\mu|$, then $|\mu| = \o f(x)$ and $x$ is a transitive
point for the restriction $f| |\mu|$.
Thus, if $Gen(\mu) \cap |\mu|$ is nonempty, e.g. if $\mu$ is ergodic, then $|\mu|$ is a topologically transitive subset
of $X$.
\end{itemize}
\end{prop}

{\bfseries Proof:}  (a) Assume  $x \in Gen(\mu)$.

If $y \not\in \o f (x)$, then choose $u : X \to [0,1]$  continuous with $u(y) = 1$ and
the closure of $\{ u > 0 \}$ disjoint from
$\o f(x)$. It follows that $u(f^i(x)) > 0$ for only finitely many $i$ and so
$\lim_{n \to \infty} \ \frac{1}{n} \Sigma_{i = 0}^{n-1} u(f^i(x)) = 0$.
Since $x \in Gen(\mu)$,  $\int u \ d\mu = 0$. Hence, $u = 0$ on the support of $\mu$. In particular, $y \not\in |\mu|$.

Because $\o f(x)$ is a chain transitive subset, it is contained in some basic set $B$.  Since distinct basic sets are disjoint,
$B = B(\mu)$ is uniquely defined by the conditions $|\mu| \subset \o f(x) \subset B(\mu)$ for $x \in Gen(\mu)$. In particular, if a basic set
$B_1$ meets $|\mu|$, then
$B(\mu) \cap B_1 \not= \emptyset$ and so $B(\mu) = B_1$. If $x_1 \in Bas(B_i) \cap Gen(\mu)$, then
 $\o f(x_1) \subset B(\mu) \cap B_1$ and so $B(\mu) = B_1$.

(b) If $x \in |\mu|$ then $\o f(x) \subset |\mu|$ because the support is closed and invariant.
The reverse inclusion follows from (a).  Thus, if
$x \in  Gen(\mu) \cap |\mu|$, then $|\mu| = \o f(x)$ and  $x$ is a transitive point for the subsystem $f||\mu|$.

$\Box$ \vspace{.5cm}

\begin{ex}\label{tracex00ac} The inclusion $|\mu| \subset \o f(x)$ for $x \in Gen(\mu) \setminus |\mu|$
may be proper. \end{ex}

{\bfseries Proof:}
If $(X,f)$ is \emph{uniquely ergodic}\index{dynamical system!uniquely ergodic}, i.e. it admits a unique invariant measure $\mu$, then the sequence
$\{ \frac{1}{n} \Sigma_{i=0}^{n-1} \d_{f^i(x)} \}$ converges to $\mu$ for every $x \in X$, i.e. $Gen(\mu) = X$.
In \cite{A93} Theorem 9.2 a non-trivial, topologically transitive system $(X,f)$ is
described with a fixed point $e$ such that
$\d_e$ is the unique invariant measure. If $x$ is a transitive point, then $x \in Gen(\d_e)$ with
$\{ e \} = |\d_e|$, but $\o f(x) = X$.

$\Box$ \vspace{.5cm}

If $f^n(x) = x$, i.e. $x$ is a periodic point, then $\frac{1}{n} \Sigma_{i= 0}^{n-1} \d_{f^i(x)}$
is the uniform
invariant measure concentrated on the periodic orbit.  It is an ergodic measure. If
an ergodic measure is not thus concentrated on
a periodic orbit, then it is
non-atomic, i.e. countable sets have measure zero.

We will call $(X,f,\mu_0)$ a \emph{measured dynamical system}\index{dynamical system!measured} when $(X,f)$ is a dynamical
system and $\mu_0$ is a \emph{background measure}\index{measure!background} on $X$.  That is, $\mu_0$  is a
full measure on $X$, so that the empty set is
the only open subset of measure zero. We also assume that $\mu_0$ is complete, i.e.
any subset of a set of measure zero is
measurable with measure zero. On the other hand, $\mu_0$ is not assumed to be $f$ invariant.

For a measured dynamical system, we call a basic set $B$ \emph{$\mu_0$-visible}\index{basic set!$\mu_0$-visible}
when $\mu_0(Bas(B)) > 0$. Since $\mu_0$ is full, a visible basic set
is $\mu_0$-visible.

We will say that a Borel set $A \subset U$ is \emph{$\mu_0$ dense}\index{set!$\mu_0$ dense} in a nonempty open set $U$ if $\mu_0(A) = \mu_0(U)$
or, equivalently, if $\mu_0(U \setminus A) = 0$.
Since the measure is full, $A$ is dense in $U$ if it is $\mu_0$ dense in $U$. If
$\mu_0(B \cap U) > 0$ and $A$ is $\mu_0$ dense in $U$, then
$B \cap A = (B \cap U) \setminus (U \setminus A)$ has positive $\mu_0$ measure.
Since the generic sets of distinct invariant measures are disjoint,
it follows if $Gen(\mu_1) \cap U$ is $\mu_0$ dense in $U$ and if $\mu_0(Gen(\mu_2) \cap U) > 0$
then $\mu_1 = \mu_2$. In particular,
$Gen(\mu_1) \cap U$ is $\mu_0$ dense in $U$ for at most one invariant measure $\mu_1$.

\begin{df}\label{tracdef01} Let $(X,f,\mu_0)$ a measured dynamical system.\index{dynamical system!measured}
 We will call  $(X,f,\mu_0)$ \emph{tractable}\index{dynamical system!tractable} when it satisfies the following conditions.

\begin{itemize}

\item[TRAC 1] There is a finite set of ergodic invariant
measures $\{ \mu_1, \dots, \mu_k \} $ with
 $\mu_0(\bigcup_{i=1}^{k} \ Gen(\mu_i)) = 1$.

 \item[TRAC 2] There are only finitely many basic sets.

\item[TRAC 3] The union in $X$ of the set of open sets
$  U $ such that $ U \cap Gen(\mu) $  is $  \mu_0 $ dense in $  U  $ for some $  \mu \in |f_*| $
 has $\mu_0$ measure one.

\end{itemize}
\end{df}
\vspace{.5cm}

\begin{theo}\label{tracprop01aa} A tractable, measured dynamical system $(X,f,\mu_0)$
satisfies the following conditions.

\begin{itemize}

\item[TRAC 1+] There is a finite set of ergodic invariant
measures $\{ \mu_1, \dots, \mu_k \} $ with $\mu_0(Gen(\mu_i)) > 0$ for $i = 1, \dots, k$ such that \\
 $\mu_0(\bigcup_{i=1}^{k} \ Gen(\mu_i)) = \sum_{i=1}^k \mu_0(Gen(\mu_i)) = 1$.
 We call $\{\mu_1, \dots, \mu_k \}$ the \emph{special ergodic measures}\index{measure!special ergodic} for the system.
For any other invariant measure $\mu$ for $f$, $\mu_0(Gen(\mu)) = 0$.

\item[TRAC 2+] There are finitely many basic sets and a basic set $B$ is  $\mu_0$-visible
iff there exists $i = 1, \dots, k$ such that $|\mu_i| \subset B$.  That is
$\{ B(\mu_i) : i = 1, \dots, k \}$ is the set of
$\mu_0$-visible basic sets. In fact,
\begin{equation}\label{eqtrac0001}
\mu_0(Bas(B)) \ = \ \sum \{ \mu_0(Gen(\mu_i)) : B = B(\mu_i) \}
\end{equation}

In particular, if $B$ is a terminal basic set, and so is visible,  then there exists $i$ such that $B = B(\mu_i)$.

\item[TRAC 3+] For $i = 1, \dots, k$ let $U_i$ be the union of those open sets $ U $ such that
$ U \cap Gen(\mu_i) $  is $  \mu_0 $ dense in $  U $.
  The sets $\{U_1, \dots, U_k \}$ are pairwise disjoint with a dense union and satisfy
  \begin{align}\label{eqtrac0001a}
  \begin{split}
  \mu_0(U_i) \ = \ \mu_0(Gen(\mu_i) &\cap U_i) \ = \ \mu_0(Gen(\mu_i)), \\
  \mu_0(Gen(\mu_j) &\cap U_i) \ = \ 0
  \end{split}
  \end{align}
  for $i,j = 1, \dots, k$ with $i \not= j$.  In addition, for $i  = 1, \dots, k$
   \begin{equation}\label{eqtrac0001b}
   U_i \ = \ X \setminus \overline{\bigcup_{j \not= i} U_j} \ = \ (X \setminus \bigcup_{j \not= i} U_j)^{\circ}.
   \end{equation}
  For an open  $U \subset X$ the following are equivalent.
  \begin{itemize}
  \item[(i)] $U$ is contained in $U_i$.
  \item[(ii)] $U \cap Gen(\mu_i)$ is $\mu_0$ dense in $U$.
  \item[(iii)] $\mu_0(U \setminus Gen(\mu_i)) = 0$.
  \item[(iv)]   $\mu_0(U \cap Gen(\mu_j)) = 0$ for all $j \not= i$.
     \end{itemize}

   We call $\{U_1, \dots, U_k \}$  the \emph{special open sets }for the system.
\end{itemize}
 \end{theo}

{\bfseries Proof:} TRAC 1+: By discarding any measures with $\mu_0(Gen(\mu)) = 0$ from
the list in TRAC 1, we may assume that
$\mu_0(Gen(\mu_i)) > 0$ for all $i$.  For any other invariant measure $\mu$, the set
$Gen(\mu)$ is disjoint from $\bigcup_i Gen(\mu_i)$
and so has $\mu_0$ measure zero.

TRAC 2+: The same basic set $B$ may contain the support of more than one special measure. That is, $B$ may equal $B(\mu_i)$ for more than one
special measure $\mu_i$.
In any case, $B = B(\mu_i)$ implies
$Gen(\mu_i) \subset Bas(B)$. Since $Gen(\mu_i)$ and $Gen(\mu_j)$ are disjoint if $i \not= j$, we clearly have
$\mu_0(Bas(B)) \ \geq \ \sum \{ \mu_0(Gen(\mu_i)) : B = B(\mu_i) \} $.  Since basins
of distinct basic sets are disjoint, the sum over all the basic sets $B$ the
left hand side of (\ref{eqtrac0001}) is at most one.  On the other hand, by TRAC 1,
the sum over all the basic sets $B$ on the right hand side equals one. This implies equality in
(\ref{eqtrac0001}) for each $B$.  In particular, if $B \not= B(\mu_i)$ for $i = 1, \dots, k$, then
$\mu_0(Bas(B)) = 0$ and $B$ is not $\mu_0$-visible.

If $B$ is terminal, then  Theorem \ref{tractheo00a} implies that $B$ is an attractor and so is visible. Since $\mu_0$ is full, it is
$\mu_0$-visible

TRAC 3+: If an invariant measure $\mu$ is not one of the special ergodic measures, then $\mu_0(Gen(\mu)) = 0$ and so
$U \cap Gen(\mu)$ cannot be $\mu_0$ dense in any nonempty open set $U$.

Since $X$ is second countable, the open set $U_i$ is a countable union
of open sets $U$ such that $U \cap Gen(\mu_i)$ is $\mu_0$ dense
in $U$. Hence, $U_i \setminus Gen(\mu_i)$ is the countable union of the sets
$U \setminus Gen(\mu_i)$ of $\mu_0$ measure zero.
Thus, $U_i$ is the maximum open set $U$ for which $U \cap Gen(\mu_i)$ is $\mu_0$ dense
in $U$. It follows that $U_i$ and $U_j$ are disjoint for $i \not= j$. Hence,
$\sum_i \mu_0(U_i) \leq 1$. On the other hand,
$\sum_i \mu_0(Gen(\mu_i)) = 1$ by TRAC 1. So equality holds for the first line of (\ref{eqtrac0001a}). In particular,
$Gen(\mu_j) \setminus U_j \supset Gen(\mu_j) \cap U_i$ has $\mu_0$ measure zero.

Since $\bigcup_{i} U_i$ has $\mu_0$ measure one and $\mu_0$ is full, the union is dense.

Now let $\tilde U_i = X \setminus \overline{\bigcup_{j \not= i} U_j}$.  Since $U_i$ is open
and is disjoint from the remaining $U_j$'s
it follows that $U_i  \subset \tilde U_i$.

If $i \not= j$ then $\tilde U_i \cap \tilde U_j = X \setminus \overline{\bigcup_{i} U_i} = \emptyset$. That is, the sets
$\tilde U_1, \dots \tilde U_k$ are pairwise disjoint.

Since the union of all of the $U_j$'s has $\mu_0$
measure one, it follows that
$$ \mu_0(\tilde U_i) = \mu_0(U_i) = \mu_0(U_i \cap Gen(\mu_i)) \leq \mu_0(\tilde U_i \cap Gen(\mu_i)))
\leq \mu_0(\tilde U_i).$$
Because equality then holds it follows that $\tilde U_i \cap Gen(\mu_i)$ is $\mu_0$ dense in
$\tilde U_i$ and so $\tilde U_i \subset U_i$.

From the definition of $U_i$ and the fact that $\mu_0(U_i \setminus Gen(\mu_i)) = 0$, it is
clear that conditions
(i) - (iii) on an open set $U$ are equivalent.  Since $\mu_0(U \cap \bigcup_{j=1}^k Gen(\mu_j)) =
\mu_0(U) > 0$, $\mu_0(U \setminus Gen(\mu_i)) = \mu_0(U \cap \bigcup_{j \not= i} Gen(\mu_j))$.
So (iii) is equivalent to (iv).

$\Box$ \vspace{.5cm}

Let $h : X_1 \to X_2$ be a continuous surjection. We say that $h$ is \emph{open at $x \in X_1$}\index{map!open at $x$}
if for every $A$ with $x \in A^{\circ}$,
$h(x) \in h(A)^{\circ}$. Observe that if $h^{-1}(y) \subset U$ for $y \in X_2$ and $U$ open
in $X_1$, then by compactness there exists
$V$ open in $X_2$ with $y \in V$ and $h^{-1}(V) \subset U$. In particular, if $h^{-1}(y)$
is the singleton $\{ x \}$, then $h$ is open
at $x$. Furthermore, the set of $y \in X_2$ such that the diameter of $h^{-1}(y)$ is less
than $\ep$ is open for every $\ep > 0$. Hence,
the set of $y$ such that $h^{-1}(y)$ is a singleton is a $G_{\delta}$ as is its preimage:
   \begin{equation}\label{eqtrac0001c}
   Inj_h \ = \ \{ x \in X_1 : h^{-1}(h(x)) = \{ x \} \}.
   \end{equation}

   \begin{prop}\label{tracprop01ab} Let $h : X_1 \to X_2$ be a continuous surjection.
   The following conditions are equivalent and
   when they hold we say that
 $h$ is \emph{almost open}\index{map!almost open}
 \begin{itemize}
 \item[(i)] If $A \subset X_1$, then $A^{\circ} \not= \emptyset$ implies $h(A)^{\circ} \not= \emptyset$.

 \item[(ii)] If $B \subset X_2$ is dense, then $h^{-1}(B)$ is dense in $X_1$.

 \item[(iii)] For every open $U \subset X_1$, $U \cap h^{-1}(h(U)^{\circ})$ is dense in $U$.

 \item[(iv)] The set $\{ x : h$ is open at $x \}$ is a dense, $G_{\d}$ subset of $X_1$.
 \end{itemize}
 \end{prop}

  {\bfseries Proof:} $(i) \Leftrightarrow (ii)$: $B$ is dense if and only if it meets every set with a
  nonempty interior and $B$ meets
  $h(A)$ if and only if $h^{-1}(B)$ meets $A$.

  $(i) \Rightarrow (iii)$: If $V$ is a nonempty open subset of $U$, then $V\cap h^{-1}(h(V)^{\circ})$
  is nonempty by (i) and it
  is contained in $V \cap (U \cap h^{-1}(h(U)^{\circ}))$.

  $(iii) \Rightarrow (iv)$: For $n = 1, 2, \dots $ let $\{V_{1,n}, \dots V_{k_n,n} \}$
  be a finite  cover of $X_1$ by open sets of
  diameter less than $1/n$. It is easy to check that the set of points at which $h$ is open is
  the $G_{\d}$ set $\bigcap_n \bigcup_i (V_{i,n} \cap h^{-1}(h(V_{i,n})^{\circ}))$.
  Then (iii)and the Baire Category Theorem implies that
  the set is dense.

 $(iv) \Rightarrow (i)$:  If $h$ is open at a point $x \in A^{\circ}$, then $h(x) \in h(A)^{\circ}.$

$\Box$ \vspace{.5cm}

Call a continuous map $h : X_1 \to X_2$ \emph{almost one-to-one}\index{map!almost one-to-one} if $Inj_h$ is dense in $X_1$.
A continuous surjection $h$ is open at every point of $Inj_h$ and so an almost one-to-one
surjection is almost open by
Proposition \ref{tracprop01ab}.
If $\mu_0$ is a full measure on $X_1$, then $h$ is called $\mu_0$ \emph{almost one-to-one}\index{map!$\mu_0$ almost one-to-one}
if $\mu_0(Inj_h) = 1$. Since
$\mu_0$ is full, $\mu_0$ almost one-to-one implies almost one-to-one.

\begin{lem}\label{traclem01ab} Let $h : X_1 \to X_2$ be an almost one-to-one surjection.
If $U \subset X_1$ is open, then
\begin{equation}\label{eqtrac0001cc}
Inj_h \cap U \ = \ Inj_h \cap h^{-1}(h(U)^{\circ}) \ = \ Inj_h \cap h^{-1}(h(U)).
\end{equation}
In particular, if
$U_1$ and $U_2$ are disjoint open subsets
of $X_1$, then $h(U_1)^{\circ}$ and $h(U_2)^{\circ}$ are disjoint open subsets of $X_2$.

If $h$ is $\mu_0$ almost one-to-one for a  full measure $\mu_0$  on $X_1$, then for any open set $U \subset X_1$,
\begin{equation}\label{eqtrac0001d}
\mu_0(U) \ = \ h_*(\mu_0)(h(U)^{\circ}) \ = \ h_*(\mu_0)(h(U))
\end{equation}

%(a) If $A_1 \cap A_2$ is nowhere dense, then $h(A_1) \cap h(A_2)$ is a nowhere dense subset of $X_2$.
%
%(b) If $h : X_1 \to X_2$ is $\mu_0$ almost one-to-one and
% $\mu_0(A_1 \cap A_2) = 0$, then $h_* \mu_0(h(A_1) \cap h(A_2)) = 0$.
\end{lem}

{\bfseries Proof:}  If $x \in Inj_h \cap U$ then since $h$ is open at points of $Inj_h$, $h(x) \in h(U)^{\circ}$.
So we have $Inj_h \cap U \subset Inj_h \cap h^{-1}(h(U)^{\circ}) \subset Inj_h \cap h^{-1}(h(U))$.
If $x \in h^{-1}(h(U))$
then there exists $x_1 \in U$ such that $h(x) = h(x_1)$.  If also $x \in Inj_h$ then
$x = x_1$ and so $x \in  Inj_h \cap U$.

 If $V$ is an open subset of $h(U_1) \cap h(U_2)$ with $U_1$ and $U_2$ disjoint, then
 (\ref{eqtrac0001cc}) implies that $h^{-1}(V)$
  is disjoint from the dense set $Inj_h$. So
$V$ must be empty.

If $\mu_0(Inj_h) = 1$, then (\ref{eqtrac0001d}) follows immediately from (\ref{eqtrac0001cc}).

%(a) Suppose that $V$ is an open subset of $h(A_1) \cap h(A_2)$.
%Because $h$ is almost one-to-one, $h^{-1}(V) \cap Inj_h$ is dense in the open set $h^{-1}(V)$. By definition of $Inj_h$,
%$h^{-1}(V) \cap Inj_h$ is a subset of the closed set $ A_1 \cap A_2$. Hence, $h^{-1}(V) \subset  A_1 \cap A_2$.
%This intersection is nowhere dense and so $h^{-1}(V) = \emptyset$. Since $V \subset h(X_1)$, $V = h(h^{-1}(V) = \emptyset$.
%
%(b) If $h(A_1) \cap h(A_2)$ had $h_* \mu_0$ positive measure then $h^{-1}(h(A_1) \cap h(A_2))$
%would have positive $\mu_0$ measure. Since $Inj_h$ has $\mu_0$ measure one, $Inj_h \cap h^{-1}(h(A_1) \cap h(A_2))$ has
%positive measure.  But $Inj_h \cap h^{-1}(h(A_1) \cap h(A_2)) \subset A_1 \cap A_2$ which has measure zero by assumption.

$\Box$ \vspace{.5cm}

If $F_1$ and $F_2$ are relations on $X_1$ and $X_2$ respectively, then we say that a continuous map
$h : X_1 \to X_2$ maps $F_1$ to $F_2$ if $h \circ F_1 \circ h^{-1}$ ($= (h \times h)(F_1) $) is
 a subset of $F_2$ or,
equivalently, $h \circ F_1 \subset F_2 \circ h$.
If $F_1$ and $F_2$ are themselves maps then this implies $h \circ F_1 = F_2 \circ h$
since inclusion implies equality for maps.

\begin{prop}\label{tracprop00} Let $(X_1,f_1)$ and $(X_2,f_2)$ be dynamical systems and let
%so that $f_1$ and $f_2$ are continuous maps.
 $h : X_1 \to X_2$ be a continuous surjection
mapping $f_1$ to $f_2$.

(a) $h$ maps $f^{-1}_1$ to $f^{-1}_2$, and for $\A = \O, \o, $ and $\CC$, $h$ maps $\A f_1$ to $\A f_2$.

(b) For $\A = \O, \o, $ and $\CC$, $h(|\A f_1|) \subset |\A f_2|$.

(c) If $B$ is an $f_1$ basic set, then there exists a unique $f_2$ basic set $h_*(B)$ such that $h(B) \subset h_*(B)$.

(d) If $B_2$ is an $f_2$ basic set, then there exists an $f_1$ basic set $B$
such that $h_*(B) = B_2$, i.e. the set map
$h_*$ on basic sets is surjective.
If $B_2$ is a terminal basic set then $B$ can be chosen to be terminal.

(e) If $\mu$ is an $f_1$ invariant measure, then $h_*\mu$ is an $f_2$ invariant measure
with $h_*\mu$ ergodic if $\mu$ is.
Furthermore, $Gen(\mu) \subset h^{-1}(Gen(h_*\mu))$. In particular, if
$Gen(\mu) \not= \emptyset$, then $Gen(h_*\mu) \not= \emptyset$
and $h(B(\mu))  \subset B(h_*(\mu))$.

\end{prop}

{\bfseries Proof:} (a) is an easy exercise and clearly implies (b).

 (c) From (b), $h(B) \subset |\CC f_2|$ and from (a) $h$ maps $\CC f_1 \cap \CC f_1^{-1}$ to
 $\CC f_2 \cap \CC f_2^{-1}$. Hence, the points of $h(B)$ lie in a single basic set $h_*(B)$.
 Since distinct basic sets are
 disjoint, $h_*(B)$ is uniquely determined by $B$.

 (d) The basic set $B_2$ is $f_2$ invariant and so $h^{-1}(B_2)$ is at least $f_1$ $^+$invariant.
 Since $h$ is surjective, $h^{-1}(B_2)$ is nonempty. Let $x \in h^{-1}(B_2)$.  Because  $h^{-1}(B_2)$ is closed and $^+$invariant
 $\o f_1 (x) \subset h^{-1}(B_2)$. By (a) $h(\o f_1 (x)) \subset \o f_2 (h(x)) \subset B_2$. There exists a basic set $B \supset \o f_1(x)$.
 By (c) the $f_2$ basic set $h_*(B)$  contains $h(B)\supset h(\o f_1(x))$ and so meets $B_2$. Again $h_*(B) = B_2$
because distinct basic sets are disjoint.

 If $B_2$ is terminal, and so is $\CC f_2 $ invariant, then $ h^{-1}(B_2)$ is $\CC f_1$ $^+$invariant and nonempty.
 So it contains a minimal nonempty $\CC f_1$ $^+$invariant set $B_1$ and this is a terminal basic set. $ h(B_1) \subset h_*(B_1) \cap B_2$
 and so $B_2 = h_*(B_1)$.

(e) Since $h \circ f_1 = f_2 \circ h$ it follows that $h_*\mu = (h \circ f_1)_*\mu = (f_2)_*h_*\mu$ if $\mu$ is
$f_1$ invariant. If $B \subset (f_2)^{-1}(B)$, then $h^{-1}(B) \subset (f_1)^{-1}(h^{-1}(B))$ and
so $\mu(h^{-1}(B))$ equals zero or one
if $\mu$ is ergodic.  Hence, $h_* \mu$ is then ergodic. Since $h_*$ is continuous and affine,
$\frac{1}{n} \sum_{i=0}^{n-1} \d_{(f_2)^i(h(x))} = h_*( \frac{1}{n} \sum_{i=0}^{n-1} \d_{(f_1)^i(x)} )$
converges to $h_*\mu$ if
$x \in Gen(\mu)$. For $B(\mu)$, the basic set which contains $|\mu|$, $h(B(\mu))$
contains $h(|\mu|) = |h_*\mu|$ and so is contained in
$B(h_*\mu)$.

$\Box$ \vspace{.5cm}

\begin{cor}\label{traccor01aaa} If $(X,f,\mu_0)$ is a tractable measured dynamical
system and $\mu_0 \gg f_*\mu_0$, then the special
open sets satisfy $f^{-1}(U_i) \subset U_i$ and $f(U_i)^{\circ} \subset U_i$ for $i = 1, \dots, k$. \end{cor}

{\bfseries Proof:}  For $j \not= i$, $U_i \cap Gen(\mu_j)$ has $\mu_0$ measure zero. By absolute continuity,
$f^{-1}(U_i) \cap f^{-1}(Gen(\mu_j)) \supset f^{-1}(U_i) \cap Gen(\mu_j)$ has $\mu_0$
measure zero. By TRAC 3+ of Theorem \ref{tracprop01aa}
it follows that $f^{-1}(U_i) \subset U_i$.

If $f(U_i)^{\circ}$ meets $\overline{\bigcup_{j \not= i} U_j}$ then it meets some
$U_j$ with $j \not= i$. This contradicts $f^{-1}(U_j) \subset U_j$.
So $f(U_i)^{\circ} \subset U_i$ follows from (\ref{eqtrac0001b}).

$\Box$

{\bfseries Remark:}  If, in addition, $f$ is an almost open surjection, then
from Proposition \ref{tracprop01ab} it follows that
$f(U_i)^{\circ}$ is dense in $f(U_i)$ from which it follows that
$f(\overline{U_i}) = \overline{f(U_i)} \subset \overline{U_i}$.

 \vspace{.5cm}

\begin{theo}\label{tractheo02} Let $(X_1,f_1)$ and $(X_2,f_2)$ be dynamical systems and
let $h : X_1 \to X_2$ be a continuous surjection
mapping $f_1$ to $f_2$.  If $(X_1,f_1,\mu_0)$ is a measured dynamical system
satisfying TRAC 1 and TRAC 2, then $(X_2,f_2,h_*\mu_0)$
is a measured dynamical system satisfying TRAC 1 and TRAC 2.
If $h$ is $\mu_0$ almost one-to-one and $(X_1,f_1,\mu_0)$ is tractable, then
$(X_2,f_2,h_*\mu_0)$ is tractable.\end{theo}

{\bfseries Proof:} Because $h$ is surjective, $h_*\mu_0$ is full when $\mu_0$ is.
Assume that TRAC 1 and TRAC 2 hold for $(X_1,f_1,\mu_0)$.

By Proposition \ref{tracprop00} (d) the surjection $h_*$ maps the set of $f_1$ basic sets
onto the set of $f_2$ basic sets.
It follows that the number of $f_2$ basic sets is less than or equal
to the number of $f_1$ basic sets and so TRAC 2 for $f_2$ follows from TRAC 2 for $f_1$.

Different special ergodic measures on $X_1$ may map to the same measure measure on $X_2$.
Let $\{\nu_1, \dots, \nu_{\ell} \}$ list the
measures in the set $\{ h_*\mu_1, \dots, h_*\mu_k \}$. By Proposition \ref{tracprop00} (e)
each $\nu_p$ is an ergodic invariant measure for
$f_2$ and $h^{-1}(Gen(\nu_p))$ contains $\bigcup \{ Gen(\mu_i) : h_*\mu_i = \nu_p \}$. Hence,
$h_*\mu_0(Gen(\nu_p)) \geq \sum \{ \mu_0(Gen(\mu_i)) : h_*\mu_i = \nu_p \}$. Summing on $p$
we see that $\sum_p h_*\mu_0(Gen(\nu_p)) = 1$ and
so
 \begin{equation}\label{eqtrac0001e}
h_*\mu_0(Gen(\nu_p)) \ = \ \sum \{ \mu_0(Gen(\mu_i)) : h_*\mu_i = \nu_p \}.
\end{equation}
Thus, TRAC 1 for $(X_2,f_2,h_*\mu_0)$ follows. Furthermore, since $\mu_0(Gen(\mu_i)) > 0$ for all $i$,
 then $h_*(Gen(\nu_p)) > 0$ for all $p$.  That is, $ \{  \nu_1, \dots, \nu_{\ell} \}$
 are the special ergodic measures for $(X_2,f_2,h_*\mu_0)$.

Now assume $h$ is $\mu_0$ almost one-to-one and that $(X_1,f_1,\mu_0)$ satisfies TRAC 3.
Assume $U \subset X_1$ is open and  $Gen(\mu) \cap U$ is
$\mu_0$ dense in $U$ for some $\mu \in |(f_1)_*|$, then let $V = h(U)^{\circ}$. From
(\ref{eqtrac0001cc}) it follows that
$Inj_h \cap U = Inj_h \cap h^{-1}(V)$. Hence,
$\mu_0(h^{-1}(Gen(h_*\mu)) \cap h^{-1}(V)) \geq \mu_0(Gen(\mu)) \cap h^{-1}(V)\cap Inj_h)$
$=  \mu_0(Gen(\mu)) \cap U \cap Inj_h)  = \mu_0(Gen(\mu)) = \mu_0(U) = \mu_0(h^{-1}(V)).$
Thus, the first inequality is an equation and $Gen(h_*\mu) \cap V$ is $h_*\mu_0$ dense in $V$.
The union of the $U$'s has $\mu_0$ measure one
by TRAC 3 and its intersection with $Inj_h$ equals that of the $h$ preimage of the corresponding $V$'s.
Hence the union of the latter
has $h_*\mu_0$ measure one.

If $V_p$ is the maximum open set such that $V_p \cap Gen(\nu_p)$ is $h_*\mu_0$ dense in $V_p$,
then applying the above arguments
with $\mu = \mu_i$ we see that
 \begin{equation}\label{eqtrac0001f}
 \bigcup \{ h(U_i)^{\circ} : h_*\mu_i = \nu_p \} \ \subset \ V_p.
 \end{equation}

$\Box$ \vspace{.5cm}

The basic set condition TRAC 2 is clearly  invariant via a conjugacy homeomorphism $h$,
but the other conditions are not except with the change in background measure from $\mu_0$ to $h_*\mu_0$.
In general, replacing $\mu_0$ by an absolutely equivalent background measure preserves tractability.
So the natural conjugacies for this situation are homeomorphisms $h$ such that $h_*\mu_0$
is absolutely equivalent to $\mu_0$.
\vspace{.5cm}

{\bfseries Example:} Suppose $X$ is a finite dimensional, connected p.l. manifold without
boundary. By the Oxtoby-Ulam Theorem,
see \cite{OU}  and \cite{AP}, any two full, non-atomic measures on $X$ are homeomorphically equivalent. Now suppose
that $(X,f)$ is a tractable system with respect to such an Oxtoby-Ulam measure $\mu_0$. Suppose that $\mu_{k+1}$ is
a non-atomic ergodic measure not included among $\mu_i$ for $i = 1, \dots, k$. Let
$\tilde \mu_0 = \frac{1}{2} \mu_0 + \frac{1}{2} \mu_{k+1}$.
This is an Oxtoby-Ulam measure and so there exists a homeomorphism $ \tilde h$ on $X$
such that $\tilde h_*(\tilde \mu_0) = \mu_0$.
Let $\tilde f = \tilde h \circ f \circ \tilde h^{-1}$ so that $\tilde h$ maps $f$ to $\tilde f$. Clearly,
$(X,\tilde f,\mu_0)$ satisfies TRAC 1 with background measure
with special ergodic measures $\tilde h_*(\mu_1), \dots, \tilde h_*(\mu_k), \tilde h_*(\mu_{k+1})$. If
$\overline{Gen(\mu_{k+1})}$ and some $\overline{Gen(\mu_{i})}$ have overlapping interiors, then TRAC 3 fails.
Finally, if $f$ admits a sequence of distinct, non-atomic, ergodic measures
$\{ \nu_j : j \in \N \}$ then $\hat \mu_0 = \frac{1}{2} \mu_0 + \Sigma_{j=1}^{\infty} \frac{1}{2^{j+1}} \nu_j$ is
an Oxtoby-Ulam measure and there exists a homeomorphism $\hat h$  mapping $\hat \mu_0$ to $\mu_0$. Let
$\hat f = \hat h \circ f \circ \hat h^{-1}$. TRAC 1 does not hold for $(X, \hat f, \mu_0)$.

$\Box$

\vspace{1cm}

\section{Subshifts from Relations and the Two Alphabet Model}
\vspace{.5cm}

For $K$  a finite set, we let $ K^{\Z_+}$ be the compact metric space of sequences $\ss = (s_0, s_1, \dots )$ in $K$.
On it the shift map $S$ is defined by $S(\ss)_i = s_{i+1}$. We define the metric $d$ by
\begin{equation}\label{alp01}
d(\ss, \tt) = inf \{ 1/(k+1) : s_i = t_i \ \text{for all } \ i < k \}.
\end{equation}
So if $s_0 \not= t_0$ then $d(\ss, \tt) = 1$.

The pair $(K^{\Z_+},S)$ is called the
\emph{full shift} with alphabet $K$.

If $G$ is a relation on a finite set $K$ we define the associated
\emph{sample path space} $K_G$ as the closed, $^+$invariant subset
\begin{equation}\label{alp02}
K_G \ = \ \{ \ss \in K^{\Z_+} : (s_i,s_{i+1}) \in G \ \text{for all} \ i \in \Z_+ \}
\end{equation}
The pair $(K_G, S)$ is a subshift of finite type.
Let $p_i : K_G \to K$ denote the projection to the $i^{th}$ coordinate,
i.e. $p_i(\ss) = s_i$. Each $p_i$ maps the shift map $S$ to the relation $G$.

Notice that for $\ss \in K_G$ every $s_i$ is in the \emph{domain} $Dom(G) = G^{-1}(K)$. It will be convenient to
assume $K = Dom(G)$. Otherwise, we replace $G$ on $K$ be its restriction to the intersection
of the non-increasing sequence $\{ G^{-n}(K) \}$.

If $G_1$ and $G_2$ are relations on finite sets $K_1$ and $K_2$, and $g : K_1 \to K_2$ maps $G_1$ to $G_2$,
 then the product of copies of $g$ restricts
to a continuous map $g : K_{1 \ G_1} \to K_{2 \ G_2}$ which maps the shift to the shift.

If $A$ is a subset of $K$,  then, using the restriction of
$G$ to $A$, i.e. the relation $G \cap (A \times A)$ on $A$, we obtain $A_G = (K_G) \cap (A^{\Z_+})$. The set $A_G$
is a closed, $^+$invariant subset of $K_G$. We will call $a_0 \dots a_n \in K^{n+1}$ an $A_G$ word if
$(a_{i-1},a_i) \in G \cap (A \times A)$ for $i = 1, \dots, n$.

For any $K_G$ word $a_0 \dots a_n$ we define the clopen \emph{ cylinder set}
$\langle a_0 \dots a_n \rangle$ $ = \{ \ss \in K_G : s_i = a_i $ for $i = 0,\dots,n \}$.

Since $\O G = \CC G$, a $G$ basic set $B$ is an $\O G \cap \O G^{-1}$ equivalence class in the set $|\O G|$ of
$G$ periodic points. The basic set $B$ is terminal when $\O G(B) \subset B$.

\begin{prop}\label{alpprop01} Let $G$ be a relation on a finite set $K$ with $K = Dom(G)$.

\begin{enumerate}

\item[(a)] If $\ss \in K_G$ then there exists a $G$ basic set $B$ and $k \in \Z_+$
such that for all $i \geq k$, $s_i \in B$ and
$S^i(\ss) \in B_G$. The basic set $B$ is called the \emph{ endset} of $\ss$, denoted $End(\ss)$. With $B = End(\ss)$,
$\o S(\ss) \subset B_G$.

\item[(b)] If $B \subset K$ is a $G$ basic set, then $B_G \subset K_G$ is an $S$ basic set and the restriction
$S|B_G$ is topologically transitive.
Furthermore, $B \mapsto B_G$ is a bijection from the set of $G$ basic sets in $K$ to the set of $S$ basic sets in $K_G$.

\item[(c)] If  $A \subset K$, then
 $\O G(A) \subset A$ if and only if   $\ss \in A_G \ \Leftrightarrow \ s_0 \in A$ for $\ss \in K_G$.  These conditions imply that
 $\CC S(A_G) \subset A_G$, $A_G$ is a clopen subset of $K_G$, and
 $\{ \ss : End(\ss) \subset A \}$ is an open subset of $K_G$.

\item[(d)] For every $s \in K$ there exists $\ss \in K_G$ with $s_0 = s$ and  $End(\ss)$ a terminal
$G$ basic set.

\item[(e)] For $B$ a $G$ basic set, the following are equivalent:
\begin{itemize}
\item[(i)] $B$ is a terminal $G$ basic set.

\item[(ii)] $B_G$ is a terminal $S$ basic set.

\item[(iii)] $B_G$ is a visible $S$ basic set.

\item[(iv)] $B_G$ is a clopen subset of $K_G$.

\item[(v)] $B_G$ is an attractor for $S$.

\item[(vi)] $\{ \ss : End(\ss) = B \}$ is an open subset of $K_G$.
\end{itemize}

 \item[(f)] The set $\{ \ss : End(\ss) $ is a terminal $G$ basic set $\}$ is a dense open subset of $K_G$.
\end{enumerate}
\end{prop}

{\bfseries Proof:} Since $Dom(G) = K$, any $s \in K$ can be extended to $\ss \in K_G$ with
$s = s_0$. If $s \in A \subset K$ with $A \subset G^{-1}(A)$, i.e. for any $s \in A, G(s) \cap A \not= \emptyset$,
then we can choose $\ss \in A_G$. In that case, $p_0(A_G) = A$. Note that if $A$ is a $G$ basic set and $s \in A$, then there exists
a $K_G$ word $a_0 \dots a_n$ with $a_0 = s = a_n$ and so all $a_i \in A$. In particular, $a_1 \in G(s)$.

(a) Since $K$ is finite, there exists $t \in K$ such that $s_i = t$ for infinitely many $i \in \Z_+$. Let $k$
be the smallest $i$ such that $s_i = t$. It follows that for all $j \geq k$  $s_j$ is
$\O G \cap \O G^{-1}$ equivalent to
$t$ and so all lie in a single equivalence class $B$. It follows that
$S^j(\ss) \in K_G \cap B^{\Z_+} = B_G$ for all $j \geq k$.
Since $B_G$ is closed and $S$ invariant it follows that $\o S(\ss) \subset B_G$.

(b) If $a_0 \dots a_n$ and $b_0 \dots b_m$ are finite $B_G$ words, then there exists a
finite $B_G$ word $c_0 \dots c_p$ with $c_0 = a_n$ and $c_p = b_0$. So
$$a_0 \dots a_n c_1 \dots c_{p-1} b_0 \dots b_m$$
 is a word in $B_G$ which can be extended to an element $\ss $ of $B_G$. We have
 $\ss \in \langle a_0 \dots a_n \rangle $ and
$S^{n+p}(\ss) \in \langle b_0 \dots b_m \rangle $. Hence, the restriction of $S$ to $B_G$ is topologically transitive.

It follows that
$B_G$ is contained in an $S$ basic set $Q$. Each projection $p_i$ maps $S$ to $G$ and so maps $Q$ into a $G$ basic set
which must be $B$ since $Q$ contains $B_G$. Thus, $ Q \subset K_G \cap B^{\Z_+} = B_G$. That is, $Q = B_G$.

Conversely, if $Q$ is an $S$ basic set and $\ss \in Q$, let $B = End(\ss)$. Suppose $s_j \in B$ for all $j \geq k$ and so
$S^j(\ss) \in B_G $ for $j \geq k$. We also have $S^j(\ss) \in Q$ for all $j$. Hence, $Q$ and the basic set $B_G$
intersect and so must be equal.

Since $p_0$ maps $B_G$ onto $B$, the map $B \mapsto B_G$ has inverse $Q \mapsto p_0(Q)$.

(c) Since any $s \in K$ can be extended to $\ss \in K_G$ with $s = s_0$, it easily follows that
$A \supset \O G(A)$ if and only if $\ss \in A_G$ whenever $s_0 \in A$. Now assume these conditions hold.

Assume $\ss \in A_G$ and $\tt \in \CC S(A_G)$. There is a finite chain $\zz^0, \dots, \zz^n$ in $K_G$ with
$\zz^0 = \ss$, $\zz^n = \tt$ and $d(S(\zz^{i-1}),\zz^i) < 1$, and so $z^{i-1}_1 = z^i_0$,
for $i = 1, \dots, n$. It follows that
the sequence $\zz = s_0,z^1_0, z^2_0, \dots, z^{n-1}_0, t_0, t_1, \dots$ lies in $K_G$ with $s_0 \in A$. Hence, it is
in $A_G$ and so $\tt = S^n(\zz) \in A_G$.

Since $A_G = \bigcup_{s \in A} \ \{ \langle  s  \rangle  \}$ it follows that $A_G$ is clopen in $K_G$.

$End(\ss) \subset A_G$ if and only if $\ss \in \bigcup_{n = 0}^{\infty} S^{-n}(A_G)$. Hence $End(\ss) \subset A_G$ is an open
condition.

% \vspace{.5cm}
%
% \begin{lem}\label{alplem02} Assume $s_0, \dots , s_n$ is a $K_G$ word with $s_0$ in a basic set $B$ and
% $s_n \not\in B$. The set $A = \O G(s_n)$ satisfies $\O G (A) \subset A$ and $A \cap B = \emptyset$. \end{lem}
%
% {\bfseries Proof:} Transitivity of $\O G$ implies $\O G (A) \subset A$. If $s \in A$, then since
% $(s_0, s_n), (s_n, s) \in \O G$, $s \in B$ would imply $s_n \in B$.
%
%$\Box$ \vspace{.5cm}
%
%(d) Since $Dom(G) = K$, we can extend
%$s$ to an  word with length greater than the cardinality of $K$ and so with repeats. Thus, we can begin
%with $s_0 = s$ and enter some basic set $B_1$.  If $B_1$ is not terminal then we can
%continue the $K_G$ word so as to eventually
%exit $B_1$. Then, continuing long enough we enter a basic set $B_2$ and again we can exit if $B_2$ is not
%terminal.  The process ends when we arrive at terminal basic set.
%Once we enter in a terminal basic set, we remain in it.
%This must happen because Lemma \ref{alplem02} implies there are no repeats among the basic sets $B_1, B_2, \dots$.
%Any continuation of the sequence to $\ss \in K_G$ has $s_0 = s$ and $End(\ss)$ terminal.

(d) Recall that for a closed relation $F$ on $X$ with $Dom(F) = X$ a terminal basic set is a minimal nonempty $\CC F$ $^+$invariant
subset and any closed $\CC F$ $^+$invariant subset contains a terminal basic set.  For $G$ on the finite set $K$ $\CC G = \O G$ and so
$\O G(s)$ contains a terminal basic set $B$. So there is a $K_G$ word $s_0, \dots, s_n$ with $s_0 = s$ and $s_n \in B$. We extend
the finite sequence to $\ss \in K_G$ with $s_0 = s$ and $s_n \in B$. Since $B$ is invariant, $s_j \in B$ for all $j \geq n$ and so
$End(\ss) = B$.

(e) From (b) $B_G$ is an $S$ basic set and so from (c) we have (i) $ \Rightarrow $ (ii), (iv), and (vi).
From (b) there are only finitely many $S$ basic sets and so by Theorem \ref{tractheo00a}
a terminal basic set is visible, i.e. (ii) $\Rightarrow$ (iii).
Because $B_G$ is a basic set, it is $S$ invariant. A clopen
 $^+$invariant set is inward and a clopen invariant set is an attractor. Hence, (iv) $\Rightarrow$ (v).

Thus, (i) implies (ii), (iii), (iv), (v) and (vi).

 Now assume (i) is false. By (d), for any $t \in B$ there exists $ \tt \in K_G$
 with $t_0 = t$ and with $End(\tt) = \hat B$ terminal and
 so not equal to  $ B$.

 Now let $\ss \in K_G$ with $End(\ss) = B$, e.g. $\ss \in B_G$. So there exists $N$ such
 that $s_n \in B$ for all $n \geq N$.
 Let $n$ be arbitrarily large with $n \geq N$. By the above remarks, there exists $\tt \in K_G$ with
 $s_n = t_0$ and $End(\tt)$ a terminal
 basic set $\hat B$. Concatenating
 we obtain $\zz \in K_G$ with $z_i = s_i $ for $i \leq n$ and with $S^{n+1}(\zz) = \tt$. Thus,
 $End(\zz) = End(\tt) = \hat B$ is terminal and not equal to $B$.

 It follows that $\{ \ss : End(\ss) = B \} \supset B_G$ has empty interior. The basic set $B_G$ is not visible and
 is not open nor is terminal, i.e. not(vi), not(iii), not(iv) and not(ii).  Since
 $\{ \ss : \o S(\ss) \subset B_G \} = \{ \ss : End(\ss) = B \}$ is not open, $B_G$ is not an attractor, i.e. not(v).

  Thus, not(i) implies
 not(ii), not(iii), not(iv), not(v) and not(vi).

(f) By (b) and (e) the set of visible $S$ basic sets is $\{ B_G : B $ a terminal $G$ basic set $ \}$ and for a terminal
basic set $B$ the basin  $Bas(B_G) = \{ \ss : End(\ss) = B \}$.  The density of the union
of the basins then follows from Theorem \ref{tractheo00a}.

$\Box$ \vspace{.5cm}

For a relation $F : I \to J $ with $I$ and $J$ finite sets and $Dom(F) = I$, a
\emph{stochastic cover}\index{stochastic cover} $\Gamma(F)$ \index{$\Gamma(G)$} is a $J \times I$ matrix such that
\begin{align}\label{alp03}
\begin{split}
\Gamma(F)_{j i} \ > \ 0 \quad &\Longleftrightarrow \quad (i,j) \in F, \\
 \text{and} \quad \Sigma_{j \in J} \ \Gamma(F)_{j i} \ &= \ 1 \quad \text{for all} \ i \in I.
 \end{split}
 \end{align}

 The  reversal of order is so that if $G : J \to K$ is another similar relation with stochastic cover
 $\Gamma(G)$ then the matrix product $\Gamma(G) \cdot \Gamma(F)$ is a stochastic cover of $G \circ F : I \to K$.

For the relation $G$ on a finite set $K$, with $Dom(G) = K$, a stochastic cover $\Gamma(G)$ induces
a Markov process with $\Gamma(G)_{s_2 s_1}$ the probability of moving from $s_1$ to $s_2$ in a single time step.
With $\Gamma(G)^n$ the $n$-fold product, $\Gamma(G)^n_{s_2 s_1}$ is the probability that the process beginning at
$s_1$, hits $s_2$ at the $n^{th}$ step (not necessarily for the first time).

An element $K$ is called \emph{transient}\index{transient} if it is not a member of some terminal $G$ basic set in $K$.
We will let $tran \subset K$\index{$trans$} denote the subset of transient elements.
If $\ss \in K_{G}$ with $s_i$ in a terminal basic set $B$ then $s_j \in B$ for all $j > i$.
Contrapositively, if $s_j$ is transient then $s_i$ is transient for all $i < j$. For any positive integer $n$
and $s \in K$, define
\begin{equation}\label{alp04}
(\Gamma(G)^n)_{tran \ s} \ = \ \Sigma \{ (\Gamma(G)^n)_{s_1 s} : s_1 \in tran \}.
\end{equation}

\begin{lem}\label{alplem03} There exists $0 < \rho < 1$ and a positive integer $n$ such that
$(\Gamma(G)^{m})_{tran \ s} \leq \rho^m $ for all $s \in K$ and positive integers $m \geq n$. \end{lem}

{\bfseries Proof:} By Proposition \ref{alpprop01} (d) for each $s$ there is an element of $K_G$ which begins at
$s$ and eventually is in some terminal $G$ basic set. So be choosing $n$ large enough, for every $s \in K$ there exists
a $K_G$ word $s_0, \dots, s_n$ with $s_0 = s$ and $s_n$ in a terminal basic set. This implies that each
$(\Gamma(G)^n)_{tran \ s} < 1$. Let $\rho_1$ be the maximum for all $s \in K$.

We first show that for every positive integer $k$, $(\Gamma(G)^{nk})_{tran \ s} \leq \rho_1^k $ for all $s \in K$.

Observe that if $s_1$ lies in some terminal basic set then $(\Gamma(G)^{i})_{tran \ s_1} = 0$
for every positive integer $i$.
Hence,
$$(\Gamma(G)^{n(k+1)})_{tran \ s} = \Sigma_{ s_1 \in tran} (\Gamma(G)^{n})_{tran \ s_1} (\Gamma(G)^{nk})_{s_1 s}
\leq \rho_1 (\Gamma(G)^{nk})_{tran \ s}$$
and so the result for all $k > 0$ follows by induction.

Let $\rho = \rho_1^{1/2n}$. For $0 \leq i < n$, $k \geq n(k+1) \cdot \frac{1}{2n} \geq (nk + i) \cdot \frac{1}{2n}$
implies that
$\rho_1^k \leq  \rho^{nk + i}.$ So
$$(\Gamma(G)^{nk+i})_{tran \ s} = \Sigma_{ s_1 \in tran} (\Gamma(G)^{i})_{tran \ s_1} (\Gamma(G)^{nk})_{s_1 s} \leq \rho_1^k \leq  \rho^{nk + i}.$$

$\Box$ \vspace{.5cm}

A \emph{distribution}\index{distribution} $v$ on $K$ is a function $v : K \to [0,1]$ with $v(s) \geq 0$ for all $s \in K$ and
$\Sigma_{s \in K} \ v(s) \ = \ 1$. The \emph{support} of $v$ is $\{ s : v(s) > 0 \}$. We call $v$ a
\emph{positive distribution}\index{distribution!positive} when
the support is all of $K$, i.e. $v(s) > 0$ for all $s \in K$.  The distribution $v$ defines a measure on $K$ by
$v(A) = \Sigma \{ v(s) : s \in A \}$ for $A \subset K$.
 We call $v$  a \emph{stationary distribution}\index{distribution!stationary}
for $\Gamma(G)$ when $\Gamma(G) \cdot v = v$, i.e. $\Sigma_{s_1 \in K} \ \Gamma(G)_{s_2 s_1} v(s_1) = v(s_2)$.

\begin{prop}\label{alpprop04} Let $\Gamma(G)$ be a stochastic cover of a relation $G$ on a finite set $K$.

(a) If $v$ is a stationary distribution for $\Gamma(G)$ then $v(s) = 0$ for all $s \in tran$.

(b) If $B$ is a terminal basic set, then there is a unique stationary distribution $v_B$ for $\Gamma(G)$ with
support equal to $B$.

(c)  If $v$ is a stationary distribution then $v = \Sigma \{ v(B) v_B : B $ a terminal $G$ basic set $ \}$.
\end{prop}

{\bfseries Proof:} (a) Since $v = \Gamma(G)^{m}\cdot v$, Lemma \ref{alplem03} defines $n$ and $\rho < 1$ such that for
every  integer $m \geq n$, $v(tran) = \Sigma_s  (\Gamma(G)^{m})_{tran \ s} v(s) \leq \rho^m$.
Hence, $v(tran) = 0$.

(b) The restriction $\Gamma(G)_{s_1 s_2}$ with $s_1, s_2 \in B$ is a an irreducible stochastic matrix when $B$ is
a terminal basic set.  By the Frobenius theory of non-negative matrices, $1$ is the dominant eigenvalue, with
one-dimensional eigenspaces, generated by positive eigenvectors. Normalizing the right eigenvector we obtain the
unique stationary distribution with support equal to $B$, see, e.g. Appendix 2 of \cite{KT}.

(c) This is clear from (a), because the restriction of $\Gamma(G)$ to the union of the terminal basic sets is
a block-diagonal stochastic matrix.

$\Box$ \vspace{.5cm}

 For the following we refer to \cite{W} Chapter 1.

 \begin{prop}\label{alpprop05} Let $\Gamma(G)$ be a stochastic cover of a relation $G$ on a finite set $K$.

 (a) For each $s \in K$ there is a measure $\mu_s$ on $K^{\Z_+}$ with support
 $\langle s \rangle = \{ \ss \in K_G : s_0 = s \}$, uniquely
 defined so that on the cylinder set $\langle s_0 \dots s_n \rangle $
\begin{equation}\label{alp05}
\mu_s\langle s_0 \dots s_n \rangle  = \Gamma(G)_{s_n s_{n-1}} \dots \Gamma(G)_{s_2 s_1} \Gamma(G)_{s_1 s}
\end{equation}
if $s = s_0$ and $= 0$ otherwise.

(b) If $B$ is a terminal basic set with associated stationary distribution $v_B$ then
$\mu_B = \Sigma_{s \in B} \ v_B(s) \mu_s$ is an ergodic invariant measure for the shift $S$ with
support $B_G$.
\end{prop}

$\Box$ \vspace{.5cm}

For any distribution $v$ on $K$ which is positive, i.e. $v(s) > 0$ for all $s \in K$, we obtain a measure
$\mu_v = \Sigma_{s \in K} v(s) \mu_s$ with full support on $K_G$.
If there are any transient elements of $K$, then $v$ is not stationary
by  Proposition \ref{alpprop04} (a) and so the measure is not invariant. However, such measures will play the
role of the required background measure.  In particular, the notion of $\mu_v$ measure zero is independent of the
choice of positive distribution $v$ since $\mu_v(A) = 0$ if and only if $\mu_s(A) = 0$ for
all $s$. That is, if $v_1$ and $v_2$ are
positive distributions on $K$, then $\mu_{v_1} \approx \mu_{v_2}$.

\begin{df}\label{alpdef06} Let $\Gamma(G)$ be a stochastic cover of a relation $G$ on a finite set $K$.
A point $\ss \in K_G$ is called $\Gamma(G)$-\emph{generic}\index{generic!$\Gamma(G)$} if $B = End(\ss)$ is a terminal $G$ basic set
and $\ss \in Gen(\mu_B)$.
\end{df}
 \vspace{.5cm}

 \begin{lem}\label{alplem06a} If $B$ is a terminal $G$ basic set and $\ss \in Gen(\mu_B)$ then $End(\ss) = B$. \end{lem}

 {\bfseries Proof:} Let $\hat B = End(\ss)$ so that $\o S(\ss) \subset \hat B_G$. If $\ss \in Gen(\mu_B)$ then
 by  Proposition \ref{tracprop00aa}(c)   $B_G = |\mu_B| \subset \o S(\ss) \subset \hat B_G$.
 As distinct basic sets are disjoint it follows that
 $B_G = \hat B_G$ and so, by Proposition \ref{alpprop01}(b), $B = \hat B = End(\ss)$.
% with $\hat B \not= B$. Then $\hat B_G$ is a basic set which is disjoint from $B_G$.
% Choose $u : K_G \to [0,1]$ continuous with $u = 1$ on $B_G$ but $u = 0$ on  $\hat B_G$.  For large enough
% $i$, $S^i(\ss) \in B_G$ and so $u(S^i(\ss)) = 0$. Hence, $\lim_{n \to \infty} \ \frac{1}{n} \Sigma_{i = 0}^{n-1} u(S^i(\ss)) = 0$.
% On the other hand, $\int u(\tt) \mu_B(d \tt) = 1$.  Hence, $\ss \not\in Gen(\mu_B)$.

 $\Box$ \vspace{.5cm}

 \begin{theo}\label{alptheo07} Let $\Gamma(G)$ be a stochastic cover of a relation $G$ on a finite set $K$ and let
 $v$ be a positive distribution on $K$. With respect to $\mu_v$ the set of points which are not $\Gamma(G)$-generic
 has measure zero.
 \end{theo}

 {\bfseries Proof:} We will call a point of $K_G$ non-generic if it is not $\Gamma(G)$-generic.

 A set has $\mu_v$ measure zero if and only if it has $\mu_s$ measure zero for all $s \in K$. Fix such an $s$.

 A point $\ss$ is certainly non-generic if $End(\ss)$ is not terminal, or, equivalently, if $s_i \in tran$ for all
 $i$. Given $s \in K$, Lemma \ref{alplem03} easily implies that $\{ \ss : s_0 = s $ and $ s_i \in tran$ for all $i \}$ has
 $\mu_s$ measure zero. Note that the set is empty if $s$ is not itself transient. It follows that the set of
 $\ss$ such that $ End(\ss)$ is not terminal has $\mu_s$ measure zero.

 If $B$ is a terminal basic set then the Birkhoff Ergodic Theorem says that $B_G \cap Gen(\mu_B)$, which is the set of
 $\Gamma(G)$-generic points in $B_G$, has $\mu_B$ measure one. Hence, the set of non-generic points in $B_G$ has $\mu_B$
 measure zero. Since $\mu_B$ is a positive mixture of the $\mu_t$'s for $t \in B$, it follows that
 the set of non-generic points in $B_G$ has $\mu_t$ measure zero for every $t \in B$.

 Consider the list of $K_G$ words $s_0 \dots s_n$ with $s_n \in B$ and $s_i \in tran$ for all $i < n$.
 The cylinder sets $\langle s_0 \dots s_n \rangle $ associated with these words provides a countable decomposition of
 $\{ \ss : End(\ss) = B \}$. If $End(\ss) = B$ then $\ss \in \langle s_0 \dots s_n \rangle $
 with $n$ the first entry time into $B$.
 To show that the non-generic points with $End(\ss) = B$ have $\mu_s$ measure zero it
 suffices to show that the intersection
 with each of these cylinder sets has measure zero. If $s_0 \not= s$ then
 $\mu_s\langle s_0 \dots s_n \rangle = 0$. So we may assume that
 $s_0 = s$. It suffices to show that the conditional probability that
 $\ss$ is non-generic given $\ss \in \langle s_0 \dots s_n \rangle $ is zero.
 Observe that $\ss$ is $\Gamma(G)$-generic iff
 $S(\ss)$ is  $\Gamma(G)$-generic. Hence, it suffices to show that zero is the conditional
 probability that $\tt = S^n(\ss)$ is non-generic.
 given $\ss \in \langle s_0 \dots s_n \rangle $, but because the process is Markov and $s_0 = s$,
 this is the $\mu_{s_n}$ measure that
 $\tt$ is non-generic. Since $s_n \in B$ and $\tt \in B_G$, the $\mu_{s_n}$ probability that $\tt$ is non-generic is zero.

 It follows that the set $\ss$ such that either $End(\ss)$ is not terminal or $End(\ss) = B$ for some terminal basic set
 $B$ but $\ss$ is not generic
has $\mu_s$ measure zero.

$\Box$ \vspace{.5cm}

The subshift associated with a relation $G$ and a stochastic cover
$\Gamma(G)$ provides the motivating example of a tractable system.

 \begin{cor}\label{alpcor08} Let $\Gamma(G)$ be a stochastic cover of a relation $G$ on a finite set $K$ and let
 $v$ be a positive distribution on $K$. The measured dynamical system $(K_G,S,\mu_v)$ is
 tractable with the set of special ergodic measures  $\{ \mu_B $ on $B_G \}$ indexed  by the terminal $G$ basic sets $B$.
 As each such $B_G$ is an attractor the sets $\{ U_B = Bas(B_G) \}$ are the special open sets.
 \end{cor}

 {\bfseries Proof:}  TRAC 1 follows from Theorem \ref{alptheo07}. TRAC 2 follows from Proposition \ref{alpprop01} (b).
For each terminal basic set $B$ the basin $Bas(B_G) = \{ \ss : End(\ss) = B \}$ is open and
contains $Gen(\mu_B)$ by Lemma \ref{alplem06a}. As $Bas(B)$ is disjoint from $Bas(B_1)$ for $B_1$ any other
terminal basic set, it follows from TRAC 1 that $\mu_v(Bas(B_G)) = \mu_v(Gen(\mu_B))$.
 This implies TRAC 3.

$\Box$ \vspace{.5cm}

Our applications will use the so-called  \emph{two alphabet model}\index{two alphabet model}
as described in Proposition 1.4 of \cite{A99}. Assume that we are given two finite sets $K^*$ and $K$ together with
a  map $\g : K^* \to K$ and a relation $J : K^* \to K$. We assume that $J$ is a \emph{surjective relation}
meaning that $p_1(J) = K^*$ and $p_2(J) = K$, or, equivalently, $Dom(J) = K^*$ and $Dom(J^{-1}) = K$.
We define the relations
$G = \g \circ J^{-1}$ on $K$ and $G^* = J^{-1} \circ \g$ on $K^*$. Observe that $Dom(G) = K$ and $Dom(G^*) = K^*$.
Since $\g \circ G^* = G \circ \g$,  $\g$ maps $G^*$ to $G$.

\begin{prop}\label{alpprop09} For finite sets $K, K^*$ assume that $\g : K^* \to K$ is a map and $J : K^* \to K$ is a
surjective relation.
Let $G = \g \circ J^{-1}$  and $G^* = J^{-1} \circ \g$.

The map $\g$ induces a bijection from the set of $G^*$ basic sets to the set of $G$ basic sets.  If $B^*$ is a $G^*$
basic set,  then $B = \g(B^*)$ is the associated
$G$ basic set. If $B$ is a $G$ basic set then $B^* = \g^{-1}(B) \cap J^{-1}(B)$
is the $G^*$ basic set to which it is associated.

If $B^*$ and $B$ are associated basic sets then the following are equivalent:
\begin{itemize}
\item[(i)] $B^*$ is terminal, i.e. $\O G^*(B^*) \subset B^*$.
\item[(ii)] $G^*(B^*) = B^*$.
\item[(iii)] $B$ is terminal, i.e. $\O G(B) \subset B$.
\item[(iv)] $G(B) = B$.
\item[(v)] $J^{-1}(B) = B^*$.
\end{itemize}
\end{prop}

{\bfseries Proof:} Observe first that
\begin{equation}\label{alp06}
\O G = \g \circ (\O G^* \cup 1_{K^*}) \circ J^{-1}, \quad \text{and}
\quad \O G^* = J^{-1} \circ (\O G \cup 1_{K}) \circ \g.
\end{equation}
Since $\g$ maps $G^*$ to $G$, $\g(|\O G^*|) \subset |\O G|$, and
for each $G^*$ basic set $B^*$ there is a unique $G$ basic set $B$ such that $\g(B^*) \subset B$.

Let $s \in |\O G|$, i.e. $(s,s) \in \O G$. By (\ref{alp06}) there
exists $s^*, t^* \in K^*$ such that $\g(s^*) = s, s \in J(t^*)$
and  $(t^*,s^*) \in \O G^* \cup 1_{K^*}$. Hence, $(s^*,s^*) \in (\O G^* \cup 1_{K^*}) \circ J^{-1} \circ \g = \O G^*$.
Thus, $s^* \in |\O G^*|$ and so $\g(|\O G^*|) = |\O G|$. Thus,
\begin{equation}\label{alp07}
|\O G| = \bigcup \{ \g(B^*) : B^* \ \text{a} \ G^* \ \text{basic set} \}.
\end{equation}

Now let $s^* \in |\O G^*|$ with $B^*$ the $G^*$ basic set containing it. By (\ref{alp06}) again, there exist
$s,t \in K$ with $\g(s^*) = s, t \in J(s^*)$ and $(s,t) \in \O G \cup 1_K$.
Clearly, $(t,s) \in \g \circ J^{-1} = G$. Hence, $s$ and $t$ are in $|\O G|$ and lie in the same basic set $B$.
It follows that $s^* \in \g^{-1}(B) \cap J^{-1}(B)$. Thus,
\begin{equation}\label{alp08}
|\O G^*| \subset \bigcup \{ \g^{-1}(B) \cap J^{-1}(B) : B \ \text{a} \ G \ \text{basic set} \}.
\end{equation}

Finally, suppose that $s^*, t^* \in \g^{-1}(B) \cap J^{-1}(B)$ for some $G$ basic set $B$. Then $s = \g(s^*) \in B$ and
there exists $t \in B$ such that $t \in J(t^*)$. Since $s,t \in B$, $(s,t) \in \O G$. So
$(s^*,t^*) \in J^{-1} \circ (\O G ) \circ \g \subset \O G^*$. Applying this first to the case $s^* = t^*$ we see that
$s^*, t^* \in |\O G^*|$. This yields equality in (\ref{alp08}). Furthermore, applying this to the pair $(s^*,t^*)$
and to the pair with
the order reversed we see that $s^*$ and $t^*$ lie in the same $G^*$ basic set. It follows that distinct basic sets of
$G^*$ are mapped by $\g$ onto distinct basic sets of $G$.

Now assume that $B^*$ and $B$ are associated basic sets.  Since $B \subset |\O G|$ it follows that for $s \in B$,
there is a sequence $s_0,\dots,s_n$ with $n \geq 1$ such that $s_0 = s_n = s$ and
$s_{i} \in G(s_{i-1})$ for $i = 1,\dots, n$.
The members of the sequence all lie in the same $\O G \cap \O G^{-1}$ equivalence class, namely $B$. In particular,
$B \subset G(B)$. Similarly, $B^* \subset G(B^*)$.
Hence, (i) $\Leftrightarrow $ (ii) and (iii) $\Leftrightarrow $ (iv).

$G^*(B^*) = J^{-1}(\g(B^*)) = J^{-1}(B)$ since $B = \g(B^*)$. So (ii) $\Leftrightarrow$ (v).

Because $\g(J^{-1}(B)) = G(B)$ and $\g(B^*) = B$, it follows that (v) $\Rightarrow$ (iv).
Conversely, $\g(J^{-1}(B)) = B$ implies that $J^{-1}(B) \subset \g^{-1}(B) \cap J^{-1}(B) = B^* \subset J^{-1}(B)$
and so $J^{-1}(B) = B^*$. Thus, (iv) $\Rightarrow$ (v).

$\Box$ \vspace{.5cm}

We will use a \emph{special two-alphabet model}\index{two alphabet model!special} where $J$ as well as $\g$ is a map. In our applications
the finite sets $K^*$ and $K$ will be identified with certain subsets of an ambient space.  Each $s^* \in K^*$ will be
contained in a unique $s = J(s^*) \in K$ with $s$ the union of the $s^*$'s
contained therein. Thus, $J$ will be a surjective map.

To define a stochastic cover, we use a map  $\nu : K^* \to [0,1]$ such that
\begin{align}\label{alp09}
\begin{split}
\nu(s^*) > 0, \quad & \text{for all} \ s^* \in K. \\
\Sigma \{ \nu(s^*) : J(s^*) = s \} \ =  & \ 1, \quad  \text{for all} \ s \in K.
\end{split}
\end{align}
Thus, $\nu$ is a positive distribution on the elements of $J^{-1}(s)$ for each $s \in K$. We will refer to such a
map $\nu$ as \emph{distribution data}\index{distribution data} for the special two alphabet model.

For the relations $G = \g \circ J^{-1}$ and $G^* = J^{-1} \circ \g$ we define from the distribution data $\nu$
the stochastic covers
\begin{align}\label{alp10}
\begin{split}
\Gamma_{s_2 s_1} \ = \quad &\Sigma \ \{ \nu(s^*) : J(s^*) = s_1, \ \text{and} \ \g(s^*) = s_2 \}, \\
\Gamma^*_{s_2^* s_1^*} \ = \quad &\begin{cases} \nu(s_2^*) \quad \text{if} \ J(s_2^*) = \g(s_1^*), \\
0 \qquad \text{otherwise}. \end{cases}
\end{split}
\end{align}

Thus, for $\Gamma$, given
 $s_1 \in K$ we choose $s^* \in J^{-1}(s_1)$ with probability $\nu(s^*)$ and then apply $\g$.
 For $\Gamma^*$, given $s_1^*$ we apply $\g$ to get $s = \g(s^*_1)$ and then
 choose $s_2^* \in J^{-1}(s)$ with probability $\nu(s_2^*)$.

We have $\Gamma_{s_2 s_1} > 0$ if and only if $(s_1, s_2) \in G$ since otherwise we are summing over the empty set with
sum equal to zero by convention. Clearly, $\Gamma^*_{s_2^* s_1^*} > 0$ if and only if $(s_1^*, s_2^*) \in G^*$.
\begin{align}\label{alp11}
\begin{split}
\Sigma_{s_2 \in K} \ \Gamma_{s_2 s_1} \ = \ \ &\Sigma \ \{ \nu(s^*) : J(s^*) = s_1 \} \ = \ 1. \\
\Sigma_{s_2^* \in K^*} \ \Gamma^*_{s_2^* s_1^*} \ = \ \ &\Sigma \ \{ \nu(s_2^*) : J(s_2^*) = \g(s_1^*) \} \ = \ 1.
\end{split}
\end{align}

We will call $\Gamma$ and $\Gamma^*$ the \emph{stochastic covers induced by} $\nu$.\index{stochastic cover!induced}

\begin{prop}\label{alpprop10} Let $\g : K^* \to K$ and $J : K^* \to K$ be a
special two alphabet model with distribution
data $\nu$. Let $\Gamma$ and $\Gamma^*$ be the stochastic covers of $G = \g \circ J^{-1}$ and $G^* = J^{-1} \circ \g$
induced by $\nu$.

If $v : K \to [0,1]$ is a stationary distribution for $\Gamma$ then $v^* : K^* \to [0,1]$ defined by
$v^*(s^*) = v(J(s^*)) \nu(s^*)$ is a stationary distribution for $\Gamma^*$.
\end{prop}

{\bfseries Proof:} Since $v$ is stationary for $\Gamma$, $\Sigma_{s_1} \ v(s_1) \Gamma_{s_2 s_1} =  v(s_2) $.

Fix
 $s_2^* \in K^*$ and let $s_2 = J(s_2^*)$
\begin{align}\label{alp12}
\begin{split}
(\Gamma^* \cdot v^*)(s_2^*) \ = \ &\Sigma_{s_1^*} \ \{\nu(s_1^*) v(J(s_1^*)) \nu(s_2^*) : J(s_2^*) = \g(s_1^*) \} \\
= \  \nu(s_2^*) \Sigma_{s_1} v(s_1)  \ &\Sigma \ \{\nu(s_1^*)   : s_2 = \g(s_1^*)  \ \text{and} \ J(s_1^*) = s_1 \} \\
 = \ \nu(s_2^*) &\Sigma_{s_1} v(s_1) \Gamma_{s_2 s_1} = \nu(s_2^*) v(s_2) = v^*(s_2^*),
 \end{split}
 \end{align}
as required.

$\Box$ \vspace{.5cm}

{\bfseries Remarks:} (a) Alternatively,  for $v^*$ and $s \in K$, the $v^*$ probability that $J(s^*) = s$, which we denote
$v^*(s)$, equals $\Sigma \ \{ v^*(s^*) : J(s^*) = s \}$ and this is $v(s)$.  So given
$J(s^*) = s$ the conditional probability of $s^*$ is $\nu(s^*)$.

(b) If $B$ and $B^*$ are associated terminal basic sets, then the positive
stationary distribution $v_B$ for $\Gamma$ with support $B$ given by Proposition \ref{alpprop04} (b)
satisfies for $s_2 \in K$:
\begin{equation}\label{alp12a}
\Sigma_{s^* \in B^* \cap \g^{-1}(s_2)} \ v_B(J(s^*)) \cdot \nu(s^*) \ = \ v_B(s_2).
\end{equation}
Note that  the sum is empty and so equals $0$ when $s_2 \not\in B$.
\vspace{.5cm}

From (\ref{alp05}) and (\ref{alp10}) we have, for $s^* \in K^*$:
\begin{equation}\label{alp13}
\mu_{s^*}\langle s^*_0 \dots s^*_n \rangle  \ = \ \nu(s^*_1) \cdot \dots \cdot \nu(s^*_n)
\end{equation}
if $ s^*_0 \dots s^*_n $ is a $K^*_{G^*}$ word with $s^*_0 = s^* $, and $ = 0$
otherwise. For $s \in K$ we define the measure on $K^*_{G^*}$
\begin{equation}\label{alp14}
\mu_s \ = \ \Sigma \ \{ \nu(s^*) \mu_{s^*} : J(s^*) = s \}.
\end{equation}
Thus,
\begin{equation}\label{alp13a}
\mu_{s}\langle s^*_0 \dots s^*_n \rangle  \ = \ \nu(s^*_0) \cdot \nu(s^*_1) \cdot \dots \cdot \nu(s^*_n)
\end{equation}
if $ s^*_0 \dots s^*_n $ is a $K^*_{G^*}$ word with $J(s^*_0) = s $, and $ = 0$ otherwise.

If $B \subset K$ and $B^* \subset K^*$ are associated terminal basic sets, let $v_B$ be the positive
stationary distribution for $\Gamma$ with support $B$ given by Proposition \ref{alpprop04} (b).
From Proposition \ref{alpprop05} (b) and Proposition \ref{alpprop10} it follows that the associated
ergodic measure on $B^*_{G^*}$ is given by
\begin{equation}\label{alp14a}
\mu_{B^*} \ = \ \Sigma_{s^* \in B^*} \ v_B(J(s^*))\nu(s^*) \mu_{s^*} \ = \ \Sigma_{s \in B} \ v_B(s) \mu_s.
\end{equation}
Thus,
\begin{equation}\label{alp15}
\mu_{B^*}\langle s^*_0 \dots s^*_n \rangle  \ = \ v_B(J(s_0^*)) \cdot \nu(s^*_0)
\cdot \nu(s^*_1) \cdot \dots \cdot \nu(s^*_n)
\end{equation}
if $ s^*_0 \dots s^*_n $ is a $K^*_{G^*}$ word with $s^*_0 \in B^*$, and $ = 0$ otherwise.

\vspace{1cm}

\section{Simplicial Dynamical Systems}
\vspace{.5cm}
An $n$ simplex $s$ in a vector space is the convex hull of a finite, nonempty, affinely
independent set $V(s) = \{ v_0, \dots, v_n \}$\index{$V(s)$},
the \emph{vertices}\index{vertex}  of $s$. A simplex $s_1$ with vertices a (proper) nonempty subset of
$V(s)$ is a (proper) \emph{face} of $s$
and we write $s_1 \leq s$ (or $s_1 < s$ for a proper face). In particular,
each vertex is a $0$ simplex face of $s$. The union of the proper faces is the the \emph{boundary} $\partial s$.
We denote by $s^{\circ}$ the open simplex $s \setminus \partial s$.  It is the interior of $s$ in the $n$ dimensional
affine subspace generated
by the vertices.

A \emph{simplicial complex}\index{simplicial complex} $L$ is a
%nonempty
finite collection of closed simplices such that $s_1 < s$ and $s \in L$ implies $s_1 \in L$ and $s_1, s_2 \in L$ with
$s_1 \cap s_2$ nonempty implies $s_1 \cap s_2$ is a
common face of $s_1$ and $s_2$ so that $V(s_1 \cap s_2) = V(s_1) \cap V(s_2)$. We let $V(L)$ denote the set $\bigcup_{s \in L}  V(s)$ of vertices of $L$. The polyhedron
associated with $L$, denoted $X(L)$\index{$X(L)$}, is the union
$\bigcup_{s \in L}  s$. A space $X$ is a \emph{polyhedron}\index{polyhedron} when there exists a simplicial complex $L$ such that
 $X(L) = X$, in which case, we say that $L$ is a triangulation of a polyhedron $X$.
Every point $x \in X$ can then be written uniquely as
\begin{equation}\label{sd02a}
x = \Sigma_{v \in V(L)} b_v(x) v \hspace{2cm}
\end{equation}
 with $b_v(x) \geq 0$ for all $v$, $ \Sigma_{v \in V(L)} b_v(x) = 1$,
and $\{ v : b_v(x) > 0 \}$ the vertices of a simplex of $L$ called the \emph{carrier}\index{carrier} of $x$. The carrier $s$ of $x$ is the
unique simplex of $L$ such that $x \in s^{\circ}$.  The numbers
$\{ b_v(x) : v \in V(K) \}$ are called the \emph{barycentric coordinates}\index{barycentric coordinates} of $x$.

We define the metric $d_L$\index{$d_L$} on $X$ as the
$\ell^1$ distance between barycentric coordinates:
\begin{equation}\label{sd02}
d_L(x,y) \ = \ \Sigma_{v \in V(L)} |b_v(x) - b_v(y)|. \hspace{2cm}
\end{equation}
So the distance between distinct vertices of $L$ is $2$ and this is the $d_L$ diameter of $X$ if $L$ is non-trivial, i.e. has more than one vertex.

In general, if we have a fixed background metric $d$ on a polyhedron $X$ then the \emph{mesh} of a triangulation
$L$ of $X$, denoted $mesh_d  L$ is the maximum of the $d$-diameters of the simplices of $L$. Of course,
if we use $d = d_L$ then
the mesh is $2$.

 A
 %nonempty
 subset $L_1$ of $L$ is a complex when $s_1 < s$ and $s \in L_1$ implies $s_1 \in L_1$
 in which case we call $L_1$ a \emph{subcomplex}\index{subcomplex} of $L$. If $X_1 \subset X(K)$ we will say that $X_1$ is
 triangulated by $L$ when $X_1 = X(L_1)$ for $L_1$ a subcomplex of $K$.

 We will identify a simplex $s$ of $L$
 with the subcomplex consisting of all of the faces of $s$.

 If  $L^*$ and $L$ are triangulations of $X$, we call $L^*$ a \emph{subdivision}\index{subdivision} of $L$ if every simplex of $L$ is
 a union of simplices of $L^*$, forming a subcomplex of $L^*$. In that case, we let $J : L^* \to L$ denote the inclusion relation:
\begin{equation}\label{sd03}
J \ = \ \{ (s^*, s) : s^* \subset s \}. \hspace{2cm}
\end{equation}
Notice that $J$ is a surjective relation, i.e. $Dom(J) = L^*$ and $Dom(J^{-1})$ $ = L$.
If no simplex of $L^*$ meets two disjoint simplices of $L$ then we call $L^*$ a \emph{proper subdivision}\index{subdivision!proper} of $L$.
Clearly, if $L^*$ is a proper subdivision of $L$ then any subdivision of $L^*$ is a proper subdivision of $L$.

A map $\g : L_1 \to L_2$ between simplicial complexes is a \emph{simplicial map}\index{map!simplicial} if $\g(V(s)) = V(\g(s))$ for all $s \in L_1$.
Clearly, the dimension of $\g(s)$ is at most that of $s$ with $dim \ \g(s) = dim \ s$ if and only if $\g$ is injective on the set of vertices
$V(s)$. When this holds for all $s \in L_1$ we call $\g$\emph{ non-degenerate}\index{map!nondegenerate}.

The associated \emph{piecewise linear}\index{map!p. l.}, or p.l., map $g : X(L_1) \to X(L_2)$ is obtained by extending linearly on each simplex the map on
the vertices.  That is,
\begin{equation}\label{sd04}
g(x) = \Sigma_{v \in V(L_1)} b_v(x) \g(v). \hspace{2cm}
\end{equation}
%If $g$ is non-degenerate then $b_{g(v)}(g(x)) = b_v(x)$ for all $v \in V(K_1)$ and $x \in |K_1|$.
In general, a map $g : X_1 \to X_2$ between polyhedra is called a p.l. map when it is associated with a simplicial
map between triangulations of $X_1$ and $X_2$. While it is obvious that the composition of simplicial maps is simplicial, it
is not so obvious that the composition of p.l. maps is p.l. It is, nonetheless true, see, e.g. \cite{JFPH} Theorem 1.11. It is also
shown there that a map $g : X_1 \to X_2$ between polyhedra is p.l. exactly when $g \subset X_1 \times X_2$ is a polyhedron.

\begin{df}\label{sddef01} A \emph{simplicial dynamical system}\index{dynamical system!simplicial} is a
simplicial map $\g : L^* \to L$ with $L^*$ a
proper subdivision of $L$.
\end{df}
\vspace{.5cm}

From such a simplicial dynamical system $\g : L^* \to L$ together with the inclusion relation $J : L^* \to L$ we obtain
a two alphabet model with relations $G = \g \circ J^{-1}$ on $L$ and $G^* = J^{-1} \circ \g$ on $L^*$.

The approximation results will require a bit more background. A \emph{derived subdivision}\index{subdivision!derived} (or first derived subdivision)
$L'$ of a complex $L$ is obtained
by choosing a vertex $v(s) \in s^{\circ}$ for each $s \in L$. The set of vertices
$ \{ v(s_0), \dots, v(s_n) \} $ spans a simplex of $L'$ if and only if they can be renumbered so that $s_0 < s_1 < \dots < s_n $ in $L$.
The derived subdivision is a proper subdivision of $L$.

If $L_1$ is a subcomplex of $L$ then $L_1'$ is a subcomplex of $L'$. Hence, if $Y$ is triangulated by $L$ then it is triangulated
by $L'$.

For $Y$ a subset triangulated by $L$ we let
\begin{align}\label{sd05}%8}
\begin{split}
N^{\circ}(Y,L) \ = \ \bigcup \{ s^{\circ} : \ s \in L \ &\text{and} \ s \cap Y \not= \emptyset \} \\
 = \ \{ x : b_v(x) > 0  \ \text{for some vertex} \ &v  \text{ of} \ L \ \text{contained in} \ Y \}.
\end{split}
\end{align}\index{$N^{\circ}(Y,L)$}

Thus,  $N^{\circ}(Y,L)$ is an open subset
of $X(L)$ which contains $Y$.

\begin{lem}\label{sdlem02} If $Y_1$ and $Y_2$ are triangulated by $L$ and $L'$ is a derived subdivision then
\begin{equation}\label{sd06}%9}
N^{\circ}(Y_1,L') \cap N^{\circ}(Y_2,L') \ = \ N^{\circ}(Y_1 \cap Y_2,L').
\end{equation}
For example, if $s_1, s_2 \in L$ then $N^{\circ}(s_1,L') \cap N^{\circ}(s_2,L') \ = $ \\ $N^{\circ}(s_1 \cap s_2,L')$.
If $s_1$ and $s_2$ are disjoint then $N^{\circ}(s_1,L')$ and $N^{\circ}(s_2,L')$ are disjoint.
\end{lem}

{\bfseries Proof:} Let the carrier of $x$ in $L'$ be the simplex spanned by $v(s_0), $ $\dots, v(s_k)$ with
$s_{i-1} < s_i$ for $i = 1, \dots k$.

Assume $x \in N^{\circ}(Y_1,L')$. Since $Y_1$ is triangulated by $L'$, $v(s_{i_1}) \in Y_1$ for some
$i_1$ between $0$ and $k$. Since $Y_1$ is triangulated by $L$, $s_{i_1} \subset Y_1$. If $x \in N^{\circ}(Y_2,L')$
as well, then, similarly,  $s_{i_2} \subset Y_2$ for some $i_2$ between $0$ and $k$. Without loss of generality we may assume
$i_1 \leq i_2$ so that $s_{i_1} \subset Y_1 \cap Y_2$ and so $x \in N^{\circ}(Y_1 \cap Y_2,L')$.

$\Box$ \vspace{.5cm}

The following is a version of the standard p.l. approximation theorem.

\begin{theo}\label{sdtheo03} Let $X$ and $Y$ be polyhedra equipped
with metrics $d_X$ and $d_Y$, and let $L$ be a triangulation of $X$.
If $f : Y \to X$ is a continuous map then there exists $\d > 0$ so that if $L_1$ is any triangulation of
$Y$ with $mesh_{d_Y} L_1 < \d$ then, for $L_1'$
any derived subdivision of $L_1$, there exists a simplicial map $\g : L_1' \to L$ with
$d_Y(f(x),g(x)) \leq 2 mesh_{d_X} L$ for all $x \in Y$, where $g$ is the p.l. map associated with $\g$. \end{theo}

{\bfseries Proof:}  $\{ N^{\circ}(s, L') :s \in L \}$ is an open cover of $X$.  Let $\ell > 0$ be a Lebesgue number for the cover
and $\d > 0$ be an $\ell$ modulus of uniform continuity for $f$.  Hence, if $mesh_{d_Y} L_1 < \d$ then
for any $t \in L_1$ there exists $s \in L$ such that $f(t) \subset N^{\circ}(s, L')$. Let $\s_f(t)$ be the
intersection of all such simplices $s$. By Lemma \ref{sdlem02} $f(t) \subset N^{\circ}(\s_f(t), L')$ and $\s_f(t) \in L$
is the smallest such simplex. This minimality implies that $t \mapsto \s_f(t)$ is incidence-preserving, i.e.
$t_1 \leq t_2$ implies $\s_f(t_1) \leq \s_f(t_2)$.

Define $\g$ by associating to the vertex $v(t) \in L_1'$ an arbitrary vertex in $\s_f(t)$. If $t_0 < t_1 \dots < t_k$
then $\s_f(t_0) \leq \s_f(t_1) \dots \leq \s_f(t_k)$ and so $\g(v(t_0)), \dots, \g(v(t_k))$ are all vertices of
$\s_f(t_k)$. Hence, each simplex of $L_1'$ is mapped to a simplex of $L$. That is, $\g$ is a simplicial map.
If the $L_1'$ carrier of $x \in Y$ is the simplex with vertices $v(t_0), \dots, v(t_k)$ then $g(x) \in \s_f(t_k)$
and $f(x) \in f(t_k) \subset N^{\circ}(\s_f(t_k),L')$. Hence, $d_X(g(x),f(x)) \leq 2 mesh_{d_X} L$.

$\Box$ \vspace{.5cm}

When $Y = X$ we can choose $L_1$ to be a subdivision of $L'$ which will ensure that it is a proper subdivision of $L$.
So we obtain the following.

\begin{cor}\label{sdcor04} Let $X$ be a polyhedron equipped with a metric $d_X$ and let $L$ be any triangulation
of $X$. If $f$ is a continuous map on $X$, then there exists a simplicial dynamical system $\g : L^* \to L$
such that $d_X(f(x),g(x)) \leq 2 mesh_{d_X} L$ for all $x \in X$, where $g$ is the p.l. map associated with $\g$.
\end{cor}

$\Box$ \vspace{.5cm}

The maps constructed above are called \emph{p.l. roundoffs}\index{map!p. l. roundoff} for $f$.

A p.l. $d$-manifold is a polyhedron $X$ such that every point of $X$ has a neighborhood which is p.l. homeomorphic
to a simplex of dimension $d$. The following result, Theorem 9.3
of \cite{A99}, allows us to refine Corollary \ref{sdcor04} to obtain
approximation by non-degenerate simplicial dynamical systems in the manifold case. The proof is rather technical and
so we will simply refer the reader to \cite{A99}.

\begin{theo}\label{sdtheo05} Let $X$ be a p.l. manifold triangulated by a complex $L_0$. Assume that $\g : L^* \to L$ is
a simplicial dynamical system with $L$ a subdivision of $L_0$. There exists a subdivision $L^*_1$ of $L^*$ and
a non-degenerate simplicial map $\g_1 : L^*_1 \to L$ such that $d_L(g(x),g_1(x)) \leq 4 mesh_{d_{L_0}} L$ for all
$x \in X$, where $g$ and $g_1$ are the p.l. maps associated with $\g$ and $\g_1$, respectively. \end{theo}

$\Box$ \vspace{.5cm}

Since any subdivision of a proper subdivision is a proper subdivision, $\g_1 : L^*_1 \to L$ is a simplicial dynamical system.

While the dynamics of simplicial dynamical systems with degeneracies is considered in \cite{A99}, we will restrict attention
here to the special case of a non-degenerate system on a p.l. $d$-manifold with $d \geq 1$. Of course, a complex of dimension
zero is just a finite set on which the dynamics is simple to describe directly.

The key result is that the local inverses of a non-degenerate simplicial dynamical system are uniformly contracting.

On the spaces $\R^n$ we use the $\ell^1$ norm and  we let $\R^n_0 = \{ a \in \R^n : \Sigma_i \ a_i = 0 \}$.

\begin{lem}\label{sdlem06} Let $(P_{ij})$ be an $(m+1)\times(d+1)$ non-negative matrix with $\Sigma_i \ P_{ij} = 1$ for all $j$.
Let $P : \R^{d+1}_0 \to \R^{m+1}_0$ be the associated linear map between the subspaces with coordinate sums equal to $0$.
If there exists $\th > 0$ such that for every pair $j_1, j_2 \in \{0,\dots, d \}$ there exists $i_0 \in \{0,\dots, m \}$ such that
$P_{i_0 j_1}, P_{i_0 j_2} \geq \th$ then the norm $|| P ||$ of the linear map is bounded by $1 - (\th/d)$. \end{lem}

{\bfseries Proof:} Assume $a \in R^{d+1}_0$ and $a \not= 0$.  Let $a_+, a_-$ be the sum of the positive and of the negative
components, respectively. Thus, $|| a || = a_+ - a_-$ and $0 = a_+ + a_-$.  That is, $a_+ = -a_- = || a ||/2$ and there
are at most $d$
nonzero components of either sign.  Hence, there exist $j_1, j_2$ such that $a_{j_1}, -a_{j_2}  \geq || a ||/2d$.
Choose $i_0$ as above for $j_1, j_2$ so that
\begin{align}\label{sd07}%10}
\begin{split}
P_{i_0 j_1} |a_{j_1}| + P_{i_0 j_2} &|a_{j_2}| - \th || a ||/d \ = \\
P_{i_0 j_1} a_{j_1} - P_{i_0 j_2} &a_{j_2} - \th || a ||/d  \ \geq \\
\max( P_{i_0 j_1} a_{j_1} + P_{i_0 j_2} a_{j_2}&, -P_{i_0 j_1} a_{j_1} - P_{i_0 j_2} a_{j_2}) \ = \\
|P_{i_0 j_1} a_{j_1}& + P_{i_0 j_2} a_{j_2}|.
\end{split}
\end{align}
Hence,
\begin{align}\label{sd08}%11}
\begin{split}
||P a || = \Sigma_i &|\Sigma_j P_{ij} a_j| \ \leq \\
\Sigma_{i,j} P_{ij} |a_j| - \th ||a||/d& \ = \ (1 - (\th/d))|| a ||.
\end{split}
\end{align}

$\Box$ \vspace{.5cm}

If $L^*$ is a subdivision of $L$ then we define
\begin{equation}\label{sd09}%12}
\th(L^*,L) \ = \ \min \{ b_v(v^*)  : v^* \in V(L^*), v \in V(L), \ \text{and} \ b_v(v^*) > 0 \}.
\end{equation}
\index{$\th(L^*,L)$}
That is, $\th(L^*,L)$ is the minimum positive $L$ barycentric coordinate of a vertex of $L^*$.

\begin{lem}\label{sdlem07} Assume $L^*$ is a proper subdivision of $L$ and that $v_1^*, v_2^* \in V(L^*)$
span a simplex in $L^*$. There exists $v_0 \in V(L)$ such that \\$b_{v_0}(v_1^*), b_{v_0}(v_2^*) \geq \th(L^*,L)$. \end{lem}

{\bfseries Proof:} If, instead, $\{ v :  b_v(v_1^*) > 0 \}$ and $\{ v :  b_v(v_2^*) > 0 \}$ are disjoint, then
the carriers in $L$ of $v_1^*$ and $v_2^*$ are disjoint and so the $1$-simplex of $L^*$ which they span meets disjoint
simplices of $L$. This means that $L^*$ is not a proper subdivision.

$\Box$ \vspace{.5cm}

\begin{prop}\label{sdprop08} Let $\g : L^* \to L$ be a non-degenerate simplicial dynamical system on a polyhedron
$X$ of positive dimension $d$. For $s^* \in L^*$
let $\bar g_{s^*}$ be the affine isomorphism from $\g(s^*) \to s^* \subset X$ which is the inverse of
the restriction of the p.l. map $g$ from $s^*$ to $\g(s^*)$. With respect to the metric $d_L$ the map
$\bar g_{s^*}$ is a contraction. For $x_1, x_2$ points in the $L$ simplex $\g(s^*)$
\begin{equation}\label{sd10}%13}
d_L(\bar g_{s^*}(x_1), \bar g_{s^*}(x_2)) \ \leq \ (1 - (\th(L^*,L)/d)) d_L(x_1,x_2).
\end{equation}
\end{prop}

{\bfseries Proof:} List the vertices of $L$ beginning with those of $\g(s^*)$ so that $V(L) = \{ v_0, \dots, v_m \}$ and
$V(\g(s^*)) = \{ v_0, \dots, v_n\}$, so that $n \leq d$.
Thus,  $\{  \bar g_{s^*}(v_0), \dots, \bar g_{s^*}(v_n) \}$ is the set of $L^*$ vertices $V(s^*)$.

Define $P_{ij} = b_{v_i}(\bar g_{s^*}(v_j))$. For $x_1, x_2$ in $\g(s^*)$ let $a_1, a_2$ be the
corresponding vectors of $\g(s^*)$ barycentric coordinates. Observe that
the remaining $L$ barycentric coordinates for $x_1$ and $x_2$ are all zero. Hence the barycentric coordinates
of $\bar g_{s^*}(x_1), \bar g_{s^*}(x_2)$ are given by $P a_1, P a_2$. Hence,
\begin{align}\label{sd11}%14}
\begin{split}
d_L(x_1, x_2) \ &= \ || a_1 - a_2 ||, \\
d_L(\bar g_{s^*}(x_1), \bar g_{s^*}(x_2)) \ = \ &|| P a_1 - P a_2 || \ = \ || P (a_1 - a_2) ||.
\end{split}
\end{align}
Since $a_1 - a_2 \in \R_0^{n + 1}$ with $n \leq d$, (\ref{sd10}) follows from Lemmas \ref{sdlem07} and \ref{sdlem06}.
Note that the vertices $\bar g_{s^*}(v_j), \ j = 0, \dots n$ are the $L^*$ vertices of $s^*$ and so any two of them
span a simplex of $L^*$.

$\Box$ \vspace{.5cm}

\begin{theo}\label{sdtheo09} Let $\g : L^* \to L$ be a non-degenerate simplicial dynamical system on a polyhedron $X$
of positive dimension $d$. Let
$g$ be the associated p.l. map on $X$. Let $J : L^* \to L$ be the inclusion relation and $G^* = J^{-1} \circ \g $ the associated
relation on the finite set $L^*$.
\begin{equation}\label{sd12}%15}
H_{\g}  \ = \ \{ (\ss^* , x) \in L^*_{G^*} \times X :  g^i(x) \in s^*_i \ \text{for} \ i \in \Z_+  \}
\end{equation}
is a continuous map from $L^*_{G^*}$ onto $X$ which maps the shift $S$ to $g$.
\end{theo}

{\bfseries Proof:} It is easy to check that $H_{\g}$ is a closed subset of $L^*_{G^*} \times X$ which is
$^+$invariant with respect to the map $S \times g$.

If $x \in X$, we choose $s^*_0 \in L^*$ so that $x \in s^*_0$. Inductively, suppose we have $s^*_0, \dots, s^*_n$
such that $(s^*_{i-1}, s^*_i) \in G^*$ and $g^i(x) \in s^*_i$ for $i = 1, \dots, n$. Since $g^{n+1}(x) \in \g(s^*_n)$
we can choose $s^*_{n+1} \subset \g(s^*_n)$ with $g^{n+1}(x) \in s^*_{n+1}$ and so $(s^*_n, s^*_{n+1}) \in J^{-1} \circ \g = G^*$.
We thus obtain, by induction, a sequence $\ss \in L^*_{G^*}$ such that $(\ss , x) \in H_{\g}$.

To complete the proof it suffices to show that $H_{\g}$ is a
map.  That is, for $\ss \in K^*_{G^*}$ there is a unique $x$ such that $(\ss , x) \in H_{\g}$.

Consider the sequence of local inverses $\bar g_{s^*_i} : \g(s^*_i) \to s^*_i$ for $i = 0, \dots n$.
Because $(s^*_i, s^*_{i+1}) \in G^*$ it follows that the image, $s^*_{i+1}$ of $\bar g_{s^*_{i+1}}$ is contained
in the domain $\g(s^*_i)$ of $\bar g_{s^*_i}$. So the sequence of images
\begin{equation}\label{sd13}%16}
\{ \bar g_{s^*_{0}} \circ \bar g_{s^*_{1}}\circ \dots \bar g_{s^*_{n - 1}} \circ \bar g_{s^*_n} (g(s^*_n)) \}
\end{equation}
is a decreasing sequence of subsets of subsets of $s^*_0$. By Proposition \ref{sdprop08} the diameter $d_K$ of the
set indexed by $n$ is at most $2 \times (1 - (\th(K^*,K)/d))^{n+1}$. Hence, the intersection is a single point $x$.
Clearly, $x$ is the unique point such that $( \ss , x) \in H_{\g}$.

$\Box$ \vspace{.5cm}

We can use $\g$ to define a sequence of successively finer subdivisions of $L$. For each $s^* \in L^*$ we can use the affine isomorphism
$\g : s^* \to \g(s^*) \in L$ to pull back the $L^*$ triangulation of $\g(s^*)$ to obtain a triangulation of $s^*$.
Doing this for every $s^*$ in $L$, we obtain a proper refinement $L^{**}$ of $L^*$ with a simplicial map $\g^* : L^{**} \to L^*$
such that $\g$ and $\g^*$ have the same p.l. map $g$ on $X$. We repeat the procedure, inductively defining a simplicial
dynamical system
$\g^{*n} : L^{*(n+1)} \to L^{*n}$ with the p.l. map $g$ for all $n$. If $s^*_0, \dots, s^*_n$ is a $L^*_{G^*}$ word
then $\bar g_{s^*_{0}} \circ \dots \circ \bar g_{s^*_{n - 1}}(s^*_n)$ is a simplex of $L^{*(n+1)}$ and every simplex
of $L^{*(n+1)}$ is obtained this way. It follows from Proposition \ref{sdprop08} that
\begin{equation}\label{sd14}
mesh_{d_L} L^{*n} \ \leq \ 2 (1 - (\th(L^*,L)/d)))^{n}. \hspace{2cm}
\end{equation}

Let $Z_0 \subset X$ denote the polyhedron triangulated by the
$d-1$ skeleton of $L^*$, i.e. the union of all of the simplices of dimension less than $d$.
A point $x$ lies in $X \setminus Z_0 \ = \ \bigcup \{ (s^*)^{\circ} \}$, with $s^*$
varying over the $d$-dimensional simplices of $L^*$,
exactly when its $K^*$ carrier
is a $d$-dimensional simplex and so the carrier is the unique simplex
of $L^*$ which contains $x$.
Since $L^*$ is a subdivision of $L$, we see that $Z_0$ contains all of the simplices of $L$ with dimension less than $d$.
Inductively, if $Z_n \subset X$ is the polyhedron triangulated by the $d-1$ skeleton of $L^{*(n+1)}$ then $Z_{n-1} \subset Z_{n}$.
We thus obtain an increasing sequence of dimension $d-1$ polyhedra in $X$. A point $x \in X \setminus Z_n$ iff
the $L^{*(n+1)}$ carrier of $x$ is $d$-dimensional, or, equivalently,
for $i = 0, \dots, n$, the $L^*$ carrier of $g^i(x)$ is $d$-dimensional. Thus, for $i = 0, \dots, n$
each $g^i(x)$ is contained in a unique simplex of $L^*$.

Because of Theorem \ref{sdtheo05}, we are primarily interested here in dynamical systems on a p.l. manifold.
It will be convenient to consider a slightly more general situation.  Call a polyhedron $X$ \emph{everywhere $d$-dimensional }
when it is a union of $d$-dimensional closed simplices.  If $L$ triangulates $X$ then every simplex of $L$ is  the
face of a $d$-dimensional simplex.

For the rest of the section we will assume that $\g : L^* \to L$ is a non-degenerate simplicial dynamical system on
an everywhere $d$-dimensional polyhedron with $d$ positive. We let $K^*$ and $K$ denote the set of $d$-dimensional simplices of
$L^*$ and $L$ respectively.  Every $d$-simplex of $L^*$ is contained in a unique $d$-simplex of $L$.  That is, the inclusion
relation restricts to a map $J$ from $K^*$ onto $K$. Because $\g$ is non-degenerate, it restricts to a map from
$K^*$ to $K$.
The relations $G^* = J^{-1} \circ \g$ and $G = \g \circ J^{-1}$ restrict to
relations  on $K^*$ and
$K$.
Notice that $K^*_{G^*}$ is a closed invariant subset of $L^*_{G^*}$ and $(K^*_{G^*},S)$ is
a subshift of finite type.

Thus, $J :\  K^* \to \ K$ and $\g :  \ K^* \to \ K$ is a special two-alphabet model with associated relations
$G^*$ and $G$.

%
%
%If $h : X_1 \to X_2$ is a continuous map then $h$ is called \emph{almost one-to-one} if
%$Inj_h = \{ x \in X_1 : h^{-1}(h(x)) $ is a singleton $ \}$ is dense in $X_1$.

\begin{theo}\label{sdtheo10} Let $\g : L^* \to L$ be a non-degenerate simplicial dynamical system on an
everywhere $d$-dimensional polyhedron $X$ with
$g$ the associated p.l. map on $X$. Let $K$ and $K^*$ be the sets of $d$-dimensional simplices of $L$ and $L^*$, respectively so that
$\g(K^*) \subset K$.
Let $J : \ K^* \to \ K$ be the inclusion map and
 $G^* = J^{-1} \circ \g $ the associated
relation on the finite set $ K^*$.
\begin{equation}\label{sd15}%19}
H_{\g}  \ = \ \{ (\ss^* , x) \in \ K^*_{G^*} \times X :  g^i(x) \in s^*_i \ \text{for} \ i \in \Z_+  \}
\end{equation}
is an almost one-to-one continuous map from $ K^*_{ G^*}$ onto $X$ which maps the shift $S$ to $g$.
\end{theo}

{\bfseries Proof:} The map labeled $H_{\g}$ here is just the restriction of the map of Theorem \ref{sdtheo09}
to the invariant subset $ K^*_{ G^*}$. For $x \in X$ we can repeat the inductive construction of the
proof of Theorem \ref{sdtheo09}, choosing
 a $d$ dimensional $s^*_{i+1} \subset \g(s^*_i)$ with $g^{i+1}(x) \in s^*_{i+1}$. The resulting sequence $\ss$
 lies in $ K^*_{G^*}$. So the restriction is still onto.

 Since $X$ is everywhere $d$-dimensional, each $d-1$ dimensional polyhedron $Z_n$ is nowhere dense in $X$ and
 so $Z_{\infty} = \bigcup_n Z_n$ is of first category in $X$.
 For $x$ in the dense set $X \setminus Z_{\infty}$, $H_{\g}^{-1}(x)$ is a singleton.

 Now let $\ss^* \in$ $ K^*_{G^*}$ and let $n$ be a positive integer. The map
 $\bar g_{s^*_{0}} \circ \bar g_{s^*_{1}}\circ \dots \bar g_{s^*_{n - 1}}$ restricts to an affine isomorphism of
 $s^*_n$ onto a $d$-dimensional simplex $t$ of $L^{*(n+1)}$ contained in $s^*_0$. Each $t \cap Z_m$ for
 $m \geq n+1$ is at most $(d -1)$-dimensional and
 so is nowhere dense in $t$. For $x \in t \setminus Z_{\infty}$,  $H_{\g}^{-1}(x)$ is a singleton  $\{ \hat \ss^* \}$
 and the first $n$ coordinates of $\hat \ss^*$ agree with those of $\ss^*$. It follows that $Inj_{H_{\g}} = $
 $\{ \ss^* : H_{\g}^{-1}(H_{\g}(\ss^*)) $ is a singleton $\}$ is dense in $K^*_{G^*}$. Thus, $H_{\g}$ is almost one-to-one
 on $K^*_{G^*}$.

$\Box$ \vspace{.5cm}

Now we turn to measures. On each $d$-dimensional simplex $s$ or $s^*$ there is a Lebesgue measure
$\l_s$ or $\l_{s^*}$ normalized by $\l_{s^*}(s^*) = 1 = \l_{s}(s)$. Since an affine mapping multiplies Lebesgue measure by a constant, the determinant, it
follows that if $s^* \in$ $ K^*$ and $ s = \g(s^*) \in$ $ K$, then $g_*(\l_{s^*}) = \l_s$, with $g$ the linear
map on $s^*$ associated with $\g$.
On the other hand, if $J(s^*) = s$, i.e. $s^* \subset s$, then $\l_s(s^*) > 0$ and the measure
$\l_{s^*}$ is just the restriction of $\l_s$ to $s^*$, normalized via division by $\l_s(s^*)$.

For the special two-alphabet model given by $\g : \ K^* \to \ K$ and $J : \ K^* \to \ K$, we define the
distribution data $\nu : \  K^* \to [0,1]$ by
\begin{equation}\label{sd16}
\nu(s^*) \ = \ \l_s(s^*) \quad \text{with} \quad s = J(s^*).
\end{equation}
Let $\Gamma$ and $\Gamma^*$ be the associated stochastic covers on $K$ and $K^*$, respectively.

From (\ref{alp13}) we obtain for each $s \in$ $K$ the measure $\mu_s$ on $ K^*_{G^*}$ with
\begin{equation}\label{sd17}
\mu_{s}\langle s^*_0 \dots s^*_n \rangle  \ = \ \nu(s^*_0) \cdot \nu(s^*_1) \cdot \dots \cdot \nu(s^*_n)
\end{equation}
if $ s^*_0 \dots s^*_n $ is a $K^*_{G^*}$ word with $J(s^*_0) = s $, and $ = 0$ otherwise. If we let
\begin{equation}\label{sd18}
\langle s \rangle = \bigcup \{ \langle s^* \rangle \subset \ K^*_{G^*} : s^* \in \  J^{-1}(s)  \},
\end{equation}
then $|\mu_s| = \langle s \rangle$.

With $v_0$ a positive distribution on $K$, let
$\mu_0 = \Sigma_{s \in  K} \ v_0(s) \mu_s$

\begin{theo}\label{sdtheo11} The surjection $H_{\g}$ of Theorem \ref{sdtheo10} is a $\mu_0$ almost one-to-one map.
If $s \in$ $ K$, then
$H_{\g}(\langle s \rangle) = s$ and $(H_{\g})_* \mu_s = \l_s$. \end{theo}

{\bfseries Proof:} It is clear that $H_{\g}$ maps $\langle s \rangle$ into $s$. If $x \in s$ then there exists $s^*_0 \in$ $K^*$
with $x \in s^*_0$ and $s^*_0 \subset s$. Beginning the inductive proof from Theorem \ref{sdtheo10} with $s^*_0$ we obtain
$\ss^* \in \langle s \rangle$ such that $H_{\g}(\ss^*) = x$.

For $s^*_0 \dots s^*_n$ a $K^*_G$ word, let $s_k = J(s^*_k) $ for $k = 0, \dots, n$ so that $s_k = \g(s^*_{k-1})$ for
$k = 1, \dots, n$, and let $s_{n+1} = \g(s^*_n)$. We assume that $s_0 = s$ so that
$\langle s^*_0 \dots s^*_n \rangle \subset \langle s \rangle$.

Define
\begin{align}\label{sd19}
\begin{split}
t_0 \ = \ s^*_0 \ = \ &\bar g_{s^*_0}(s_1), \\
t_k  =   \bar g_{s^*_0} \circ \dots \circ \bar g_{s^*_{k-1}}(s^*_k) &= \bar g_{s^*_0} \circ \dots \circ \bar g_{s^*_k}(s_{k+1}),
\end{split}
\end{align}
for $k = 1, \dots, n$.

We prove by induction that $\l_s(t_k) = \mu_s( \langle s^*_0 \dots s^*_k \rangle)$.

For $k = 0$ this is the definition of $\nu(s^*_0)$.

For each $k \geq 1$, $\bar g_{s^*_0} \circ \dots \circ \bar g_{s^*_{k-1}}$ maps $\l_{s_k}$ to $\l_{t_{k-1}}$
and by induction hypothesis $\l_s(t_{k-1}) = \nu(s^*_0) \dots \nu(s^*_{k-1})$. The measure
$\l_{t_{k-1}}  $ is the restriction of  $\l_s$ divided by $\l_s(t_{k-1})$. By definition,
$\l_{s_k}(s^*_k) = \nu(s^*_k)$. That is, the $s_k$ measure of $s^*_k$ is $\nu(s^*_k)$. So the
$t_{k-1}$ measure of $\bar g_{s^*_0} \circ \dots \circ \bar g_{s^*_{k-1}}(s^*_k) = t_k$ is $\nu(s^*_k)$.

It follows that
$$\l_s(t_k) = \l_s(t_{k-1}) \cdot \nu(s^*_k) = \nu(s^*_0) \cdot \dots \cdot \nu(s^*_k) = \mu_s(\langle s^*_0 \dots s^*_k \rangle).$$

As is shown in the proof of
Theorem \ref{sdtheo10} $H_{\g}$ is injective on the points of $K^*_G \setminus (H_{\g})^{-1}(Z_{\infty})$,
and the map is open at these points. So $H_{\g}$ restricts to a homeomorphism from $ K^*_G \setminus (H_{\g})^{-1}(Z_{\infty})$
onto $X \setminus Z_{\infty}$. $\l_s(s \cap Z_{\infty}) = 0$
because $Z_{\infty}$ is a countable union of lower dimensional polyhedra.
Since the cylinders $\langle s^*_0 \dots s^*_n \rangle$ comprise a basis for $ K^*_G$, once we show that $\mu_s( (H_{\g})^{-1}(Z_{\infty})) = 0$
for all $s \in K$ it will follow that $(H_{\g})_* \mu_s = \l_s$.

%
% it will suffice to prove that
%$\mu_s( (H_{\g})^{-1}(Z_{\infty}) \cap \langle s^*_0 \dots s^*_n \rangle ) = 0$.

It suffices to show that for each $d$-simplex $t \in L^{*n} \cap s$ with $s \in K$, the boundary $ \partial t$
lifts via $H_{\g}$ to a set of $\mu_s$ measure zero.
If $x \in \partial t$ and $\ss^* \in \langle s \rangle$ with $x = H_{\g}(\ss^*)$, then for $k > n$, $t_k$ as defined above is
a simplex of $L^{*k}$ which contains $x$ and so meets the boundary $\partial t$.

Since $\l_s(\partial t) = 0$ there exists for every $\ep > 0$ a $\d > 0$ so that the $\l_s$ measure of the $\d$ neighborhood
of $\partial t$ is less than $\ep$. For sufficiently large $k$, the mesh of $L^{*k}$ is less than $\d$. So every $\ss^*$
with $H_{\g}(\ss^*) \in \partial t$ has initial word $s^*_0 \dots s^*_k$ whose cylinder is mapped by $H_{\g}$
onto such a simplex $t_k$.
The sum of the $\l_{t}$ measures of all such $t_k$'s is less than $\ep$ and it is equal to the $\mu_s$ measure of the sum of the
cylinders.  It follows that $\mu_s((H_{\g})^{-1}(\partial t)) = 0$.

$\Box$ \vspace{.5cm}

Now let $B$ and $B^*$ be associated terminal basic sets for $( K,G)$ and $( K^*,G^*)$. Let $\bar B^*$ and $\bar B$
consist of the subcomplexes of $L^*$ and $L$ consisting of all of the faces of the simplices of $B^*$ and $B$, respectively.
Because these are terminal sets, Proposition \ref{alpprop09} implies that $J^{-1}(B) = B^*$ and so $\bar B^*$ is a subdivision of $\bar B$.
They are everywhere $d$-dimensional with
 $B $ and $B^*  $ the sets of $d$ simplices of$  (\bar B)$ and $ (\bar B^*)$, respectively. Thus,
\begin{equation}\label{sd21}
X(\bar B) \ = \ \bigcup \{ s \in B \} = X( \bar B^*) \ = \ \bigcup \{ s^* \in B^* \}.
\end{equation}

There is a unique positive $\Gamma$ invariant
distribution $v_B$ on $B$, i.e. for $s \in $ $ K$, (\ref{alp12a}) becomes
\begin{equation}\label{sd20}
\Sigma_{s^* \in B^* \cap \g^{-1}(s)} \ v_B(J(s^*)) \cdot \nu(s^*) \ = \ v_B(s).
\end{equation}

Since the restriction $\g : \bar B^* \to \bar B$ is a non-degenerate
simplicial dynamical system, Theorem \ref{sdtheo10} implies that
$H_{\g}(B^*_{G^*}) = X( \bar B^*) = X(\bar B)$.

\begin{theo}\label{sdtheo12}  Let $\g : L^* \to L$ be a non-degenerate simplicial dynamical system on an
everywhere $d$-dimensional polyhedron $X$ with
$g$ be the associated p.l. map on $X$.  Let $K$ and $K^*$ be the sets of $d$-dimensional
simplices of $L$ and $L^*$, respectively so that
$\g(K^*) \subset K$. Let $J : \  K^* \to \  K$ be the inclusion map and
$G = \g \circ J^{-1}$ and $G^* = J^{-1} \circ \g $ the associated
relations on the finite sets $ K, K^*$.  With Lebesgue distribution data\index{distribution data!Lebesgue} $\nu(s^*) = \l_s(s^*)$ for $s = J(s^*)$, let
$\Gamma, \Gamma^*$ be the associated stochastic covers for $G$ and $G^*$.
With $v_0$ a positive distribution on $ K$, let
$\l_0 = \Sigma_{s \in   K} \ v_0(s) \l_s$ be the full, locally Lebesgue background measure on $X$.

If  $B$ and $B^*$ are associated terminal basic sets for $(K,G)$ and $(K^*,G^*)$, then
$\l_B = \Sigma_{s \in B} \ v_B(s) \l_s$ is an ergodic measure for $g$ with support $H_{\g}(B^*_G) = X( \bar B^*) = X(\bar B) $, where
$v_B$ is the $\Gamma$ invariant distribution on $B$.

With the measures $\{ \l_B \}$ indexed by the terminal basic sets of $(K,G)$, the system $(X,g,\l_0)$ is
tractable. The special open sets $\{ U_B \}$ are given by $U_B = H_{\g}(Bas_B)^{\circ}$ with $Bas(B)$ the basin
 of the attractor $B^*_{G^*}$ for $(K^*_{G^*},S)$.
\end{theo}

{\bfseries Proof:} By Theorem \ref{sdtheo11}, $H_{\g}$ is a $\mu_0$ almost one-to-one surjection. From Corollary
\ref{alpcor08} it follows that $(  K^*_{G^*},S)$ is tractable and from Theorem \ref{tractheo02} it then follows
that the image system under $H_{\g}$ is tractable.  From Theorem \ref{sdtheo11} it follows
that
$$(H_{\g})_*( \Sigma_{s \in B} v_B(s) \mu_s) = \Sigma_{s \in B} v_B(s) \l_s = \l_B.$$

$\Box$ \vspace{.5cm}

All of the locally Lebesgue measures obtained by taking different positive distributions $v_0$ on $K$ are
absolutely equivalent to one another and so any one can be used as the background measure. Clearly, for each
terminal $G$ basic set $B$ the measure
$\l_B$ is absolutely continuous with respect to $\l_0$. In fact, it is a locally Lesbesgue measure on the
the $d$-dimensional polyhedron $X(\bar B) = X(\bar B^*)$.

In can happen that for distinct terminal basic sets $B_1$, $B_2$ the polyhedra $X(\bar B_1)$ and $X(\bar B_2)$ have
a non-empty intersection, but only
at boundary points, i.e. the intersection is a polyhedron of dimension less than $d$. In that case, and in other cases
as well,  $X(\bar B_1)$ and $X(\bar B_2)$ are contained in the same basic set for $g$. That is, the map which takes
the $G^*$ basic set $B^*$ in $ K^*$ to the $g$ basic set which contains $X(\bar B^*)$ need not be injective.
By Proposition \ref{alpprop01} every
terminal basic set of $ K^*$ is equal to $B^*_{G^*}$ with $B^*$ some terminal $G^*$ basic set
and this has image $X(\bar B^*) = X( \bar B )$. From Proposition \ref{tracprop00} again
it follows that every terminal $g$ basic set contains  $X( \bar B )$ for some terminal $G$ basic set $B$.
A $g$ basic set which contains  $X( \bar B )$ for a  terminal $G$ basic set $X( \bar B )$ is necessarily visible. This was
part of the proof for Theorem \ref{tractheo02} of the TRAC 3 property.
However, such a $g$ basic set need not be itself terminal.
\vspace{.5cm}

{\bfseries Examples} (a) Let $X$ be the interval $[0,2]$ triangulated by $L$ with vertices $\{0, 1, 2 \}$.
Let $L'$ be the derived subdivision with additional vertices at $ \frac{1}{2}$ and $\frac{3}{2}$.
Define $\g : L' \to L$ by
\begin{equation}\label{sd22}
0,1,2 \mapsto 1, \quad \text{and} \quad \frac{1}{2} \mapsto 0, \frac{3}{2} \mapsto 2.
\end{equation}
On $K = \{ [0,1], [1,2] \}$, $G = ([0,1],[0,1]), ([1,2],[1,2]) = 1_K$. The map $g$ is chain transitive and so
$X$ is the unique $g$ basic set, which is, of course, terminal. $B_1 = \{ [0,1] \}, B_2 = \{ [1,2] \}$ are
distinct terminal $G$ basic sets with associated Lebesgue ergodic measures on each interval.

(b)  Let $X$ be the interval $[0,3]$ triangulated by $L$ with vertices $\{0, 1, 2, 3 \}$.
Let $L'$ be the derived subdivision with additional vertices at $ \frac{1}{2}, \frac{3}{2},  \frac{5}{2}$.
 Define $\g : L' \to L$ by
\begin{align}\label{sd23}
\begin{split}
0,1 \mapsto 1, \quad &\text{and} \quad \frac{1}{2} \mapsto 0, \frac{3}{2} \mapsto 2, \\
2,3 \mapsto 3, \quad &\text{and} \quad \frac{5}{2} \mapsto 2.
\end{split}
\end{align}
On $K$, $G = ([0,1],[0,1]), ([1,2],[1,2]), ([1,2],[2,3]), ([2,3],[2,3])$. The terminal $G$ basic sets are
$B_1 = \{ [0,1] \}$, and $ B_2 = \{ [2,3] \}$. The $g$ basic sets are $[0,1]$ and $[2,3]$.
The latter is terminal. The $g$ basic set $[0,1]$ is visible but not terminal.
\vspace{.5cm}

If $f$ itself is a non-degenerate p.l. map we can use a special approximation.

\begin{lem}\label{sdlem13} Let $L$ be a subcomplex of $K$. If $L_1$ is a subdivision of $L$ then there exists
a subdivision $K_1$ of $K$ which agrees with $L_1$ on $|L|$. \end{lem}

{\bfseries Proof:} List the simplices $\{ s_1, s_2, \dots s_m \}$ of $K \setminus (L \cup V(K))$ so that $dim \ s_i \ \leq \ dim \ s_j$ when
$i < j$.  Let $L'_1 = L_1 \cup V(K)$. For $n \geq 2$ we inductively define $L'_n$ to be a complex triangulating
$X(L'_1) \cup s_1 \cup \dots s_{n-1}$ with $L'_{n-1}$ a subcomplex of $L'_n$. Successively star at each simplex on the list,
introducing a vertex $v(s_n) \in s_n^{\circ}$ and using the cone
triangulation  $v(s_n) * L'_n \cap \partial s_i$ to extend $L'_n$ over $s_n$. Let $K_1 = L'_{m+1}$.

$\Box$ \vspace{.5cm}

\begin{prop}\label{sdprop14} Let $L_1$ be a subdivision of a complex $L$.  There exists a subdivision $L_2$ of $L_1$
and a non-degenerate simplicial map $\xi : L_2 \to L$ which is subordinate to the identity on $L$.  That is, with $h$
the p.l. map associated with $\xi$, $h(s) = s$ for each simplex $s \in L$. \end{prop}

{\bfseries Proof:} We proceed by induction on the dimension of $L$.

If the dimension is zero then $L_1 = L = V(L)$ and we use $L_2 = L_1$ and $\xi = 1_{V(L)}$.

Notice that if $L$ is a one-simplex then there is a non-degenerate simplicial map from $L_1$ to $L$ which fixes the endpoints
iff there are an even number of vertices. So if $L_1$ has an odd number of vertices we obtain $L_2$ by introducing one
additional vertex.

Now assume that the dimension of $L$ is $n \geq 1$. Let $\tilde L$ be the $n - 1$ skeleton and $\tilde L_1$ the
subdivision of $\tilde L$ which is the restriction of $L_1$.  By induction hypothesis there is a subdivision
$\tilde L_2$ of $\tilde L_1$ and a non-degenerate simplicial map $\tilde \xi : \tilde L_2 \to \tilde L$ which
is subordinate to the identity.

For each simplex $s \in$ $ ^nL$, $\tilde \xi$ maps $\tilde L_2|\partial s \to \tilde L|\partial s$.
$\tilde L_2|\partial s$ is a subdivision of $\tilde L_1|\partial s$ and by Lemma \ref{sdlem13} we can
extend the subdivision to a subdivision  $\hat L_2|s$ of $\tilde L_1|\partial s$. Extend $\tilde \xi$
to $\hat \xi$ from $\hat L_2$ to $L$ by mapping each vertex of $\hat L_2|s \setminus \hat L_2|\partial s$
to a vertex of $s$. The simplicial map $\hat \xi$ is non-degenerate on the $\partial s$. Apply Proposition 9.2 of
\cite{A99} to obtain a subdivision $L_2|s$ of $\hat L_2|s$ which agrees with  $\hat L_2|\partial s$ on $\partial s$
and non-degenerate simplicial
$\xi :L_2|s \to L|s$ which agrees with $\hat \xi$ on $\hat L_2|\partial s$.  This defines $\xi : L_2 \to L$ subordinate to
the identity on $L$.

$\Box$ \vspace{.5cm}

\begin{theo}\label{sdtheo15} Let $L$ and $K$ be triangulations of a polyhedron $X$ with $L_1$ a common subdivision.
If $\phi : K \to L$ is a non-degenerate simplicial map, then there exists $L^*$ a proper subdivision of $L_1$ and
a non-degenerate simplicial map $\xi : L^* \to K$ which is subordinate to the identity on $K$. The composition
$\g = \phi \circ \xi : L^* \to L$ is a non-degenerate simplicial dynamical system.\end{theo}

{\bfseries Proof:} The derived subdivision $L_1'$ is a proper subdivision of $L_1$. By Proposition \ref{sdprop14}
 there exists $L^*$ a  subdivision of $L_1'$ and
a non-degenerate simplicial map $\xi : L^* \to K$ which is subordinate to the identity on $K$. Since $L^*$ is a
subdivision of $L_1'$ it is a proper subdivision of $L_1$ and so is a proper subdivision of $L$. It follows that
$\g$ is a simplicial dynamical system.  As the composition of non-degenerate simplicial maps, it is non-degenerate.

$\Box$

\vspace{1cm}

\section{Shift-like Dynamical Systems}
\vspace{.5cm}

With $A$ a finite set, which we regard as an alphabet, let $X = A^{\Z_+}$ be the Cantor set of infinite sequences in $A$.
On it is defined the surjective shift map $\hat S$ by $\hat S(x)_i = x_{i+1}$ for all $i \in \Z_+$. We use the notation
$\hat S$ because we will be using another shift map below.

For $[n] = \{ 0,...,n-1 \}, A^{[n]}$ is the set of words of length $n$ for any  $n \in \Z_+$.
We write $n(w)$ for the length of a word $w$.
% with $n(x) = \infty$ for all $x \in X$.
Define the projection $J_n : X \to A^{[n]}$ by $J_n(x)_i = x_i$ for $i \in [n]$.  That is, $J_n(x)$ is the initial word of
length $n$ in $x$. Similarly, if $m \geq n$ we define $J_n : A^{[m]} \to A^{[n]}$.

For a word $a$ of length $n$ we define the clopen cylinder set $< a > \subset X$ by
\begin{equation}\label{sft00aa}
\langle a \rangle \ = \ \{ x \in X : J_n(x) = a \}.
\end{equation}

As before, see (\ref{alp01}), we
use the metric $d$ given by
\begin{equation}\label{sft00}
d(x, y) = inf \{ 1/(k + 1) : x_i = y_i \ \text{for all } \ i < k \}.
\end{equation}
Thus, $J_n(x) = J_n(y)$ if and only if $d(x,y) \leq 1/(n + 1)$.

We concatenate words in the obvious way.  If $s$ is a word and $x \in X$ then $sx \in X$ with
\begin{equation}\label{sft01}
(sx)_i \ = \ \begin{cases} s_i \ \ \text{for} \ i < n(s), \\ x_{i - n(s)}  \ \text{for} \ i \geq n(s). \end{cases}
\end{equation}
If $t$ is another word, we similarly define $st$ to obtain a word of length $n(s) + n(t)$.

\begin{df}\label{sftdef01} Given  positive integers $n, k$, we say that a continuous map
$f$ on $X$ is $(n,n+k)$ \emph{continuous}\index{map!$(n,n+k)$ continuous} if  $J_n(f(x))$ depends only on $J_{n+k}(x)$, i.e.
the preimage $f^{-1}(\langle J_{n}(f(x)) \rangle ) $ contains $\langle J_{n+k}(x) \rangle $.  Equivalently, the
relation $\g = J_n \circ f \circ (J_{n+k})^{-1} : A^{[n+k]} \to A^{[n]}$ is a mapping. In that case, we say that $\g$ is
\emph{associated with} $f$.

Given a map $\g : A^{[n+k]} \to A^{[n]}$, we define
the associated \emph{shift-like} map\index{map!shift-like} $g$ on $X$ by $g(s^*x) = \g(s^*)x$ for $s^* \in A^{[n+k]}$. Thus, $g$ is defined by
\begin{equation}\label{sft02}
\hat S^{n} \circ g \ = \ \hat S^{n+k}, \quad \text{and} \quad J_n \circ g = \g \circ J_{n+k},\\
\end{equation}
\end{df}
\vspace{.5cm}

\begin{prop}\label{sftprop02}
If $f$ is a continuous map on $X$ then for any positive integer $n$ there exists a positive integer $k$ so that
$f$ is $(n,n+k)$ continuous. \end{prop}

{\bfseries Proof:} The continuous map $f$ is uniformly continuous and so there exists $k$ so that
$d(x,y) \leq \frac{1}{n+k+1}$ implies $d(f(x),f(y)) \leq \frac{1}{n+1}$. From the definition (\ref{sft00}) of the metric
this says that $J_{n+k}(x) = J_{n+k}(y)$ implies $J_n(f(x)) = J_n(f(y))$.

$\Box$ \vspace{.5cm}

{\bfseries Remark:} It can happen that $f$ is $(n,n)$ continuous, i.e. we could use $k = 0$ for $n$.  For example, this is
true if $f$ is the identity map.  However, we will always choose $k > 0$.
 \vspace{.5cm}

For fixed positive integers $n, k$, let $K = A^{[n]}$, $K^* = A^{[n+k]}$ and $\bar K = A^{[k]}$
be the sets of words of length $n, n+k$
and $k$, respectively. Let $J = J_n : K^* \to K$  so that $s = J(s^*)$ is the initial word of
length $n$ and $\bar J : K^* \to \bar K$ so that  $\bar s = \bar J(s^*)$ is the terminal word of length $k$.
That is, $s^* = s \bar s$.

The set of cylinders
$\A(K) = \{ \langle s \rangle : s \in K \}$ is a clopen partition of $X$ with $y \in \langle s \rangle$ if and only if $y = sx$ with
$x = \hat S^n(y)$. The clopen partition $\A(K^*) = \{ \langle s^* \rangle : s^* \in K^* \}$ is a refinement of $\A(K)$ with
$ \langle s^* \rangle \subset \langle s \rangle $ for $(s^*,s) \in K^* \times K$ if and only if $s = J(s^*)$. That is, we can regard
$J$ as the inclusion map from $\A(K^*)$ to $\A(K)$.

Given a map $\g : K^* \to K$ we obtain a special two-alphabet model with relations $G = \g \circ J^{-1}$ and
$G^* = J^{-1} \circ \g$  on $K$ and $K^*$, respectively. We can identify $(K^*)^{\Z_+}$ with $X = A^{\Z_+}$ via
the homeomorphism which associates to $\ss^* \in (K^*)^{\Z_+}$  the obvious infinite concatenation
$s^*_0 s^*_1 \dots $ of successive words of length $n + k$. This identifies the shift $S$ on $(K^*)^{\Z_+}$
with $\hat S^{n+k}$ on $X$. As before, we restrict the shift $S$ to the closed invariant subspace
$K^*_{G^*}$.

We define the map $H_{\g} : K^*_{G^*} \to X$ as follows.\index{$H_g$}

For $\ss^* = s^*_0 s^*_1 \dots \ \in K^*_{G^*}$, let $s_j = J(s^*_j)$ and $\bar s_j = \bar J(s^*_j)$.
Thus, $s^*_j = s_j \bar s_j$ and $\g(s^*_j) = s_{j+1}$ for all $j$.  Define
\begin{equation}\label{sft06}
H_{\g}(\ss^*) \ = \ s_0 \bar s_0 \bar s_1 \bar s_2 \dots \ = \ s^*_0 \bar s_1 \bar s_2 \dots .
\end{equation}
It follows that $g(H_{\g}(\ss^*)) = s_1 \bar s_1 \bar s_2 \dots = H_{\g}(S(\ss^*))$. That is, $H_{\g}$ maps $S$ to $g$.
Furthermore,
\begin{equation}\label{sft07}
J_{n+k}(H_{\g}(\ss^*)) \ = \ s^*_0 \ = \ p_0(\ss^*),
\end{equation}
with $p_0 : (K^*)^{\Z_+} \to K^*$ the first coordinate projection.

\begin{theo}\label{sfttheo03} With $K = A^{[n]}$ and $K^* = A^{[n+k]}$, let $\g : K^* \to K$ be a map.

\begin{itemize}
\item[(a)] Assume $f$ is a map on $X$ which is $(n,n+k)$ continuous and which is associated with $\g$. Let
$R^f : X \to K^*_{G^*}$\index{$R^f$} be defined by $R^f(x)_j = J_{n+k}(f^j(x))$. The map $R^f$ is continuous and maps
$f$ on $X$ to the shift $S$ on $K^*_{G^*}$.

\item[(b)] If $g$ on $X$ is the shift-like map associated to $\g$, then $g$ is $(n,n+k)$ continuous with
$\g = J_n \circ g \circ (J_{n+k})^{-1}$, i.e. $\g$ is associated with $g$. The map $R^g$ is a homeomorphism
with inverse $H_{\g}$ and so is a conjugacy from $g$ on $X$ to $S$ on $K^*_{G^*}$.

\item[(c)] If $f_1$ is a continuous map on $X$ which is also associated with $\g$ then
\begin{equation}\label{sft03}
\begin{split}
J_n \circ f_1 \ = \ \g \circ J_{n+k} \ = \ J_n \circ f. \\
\text{In particular,} \quad J_n \circ f = J_n \circ g,
\end{split}
\end{equation}
 and so $d(f(x),g(x)) \leq \frac{1}{n+1}$ for all $x \in X$.

\item[(d)] The map $Q^f = (R^g)^{-1} \circ R^f = H_{\g} \circ R^f$\index{$Q^f$} on $X$ maps $f$ to $g$. If $y = Q^f(x)$ then for every
$j \in \Z_+$,
\begin{equation}\label{sft04}
J_{n+k}(f^j(x)) \ = \ J_{n+k}(g^j(y)),
\end{equation}
and so $d(f^j(x),g^j(y)) \leq \frac{1}{n+k+1}$ for all $j$.
\end{itemize}
\end{theo}

{\bfseries Proof:} (a) Clearly, $J \circ J_{n+k} = J_n  : X \to K$  and so for any $x$ and $j$,
\begin{equation}\label{sft05}
J(R^f(x)_{j+1}) = J_n (f^{j+1}(x)) = \g(J_{n+k}(f^j(x))) = \g(R^f(x)_j).
\end{equation}
That is, $(R^f(x)_j,R^f(x)_{j+1}) \in G^*$ and so $R^f$ maps $X$ into $K^*_{G^*}$. Since $(K^*)^{\Z_+}$ is a product of
finite discrete spaces, $R^f$ is clearly continuous and clearly maps $f$ to $S$.

(b) Since $J_n \circ g = \g \circ J_{n+k}$ and $J_{n+k}$ is onto, it follows that $\g = J_n \circ g \circ (J_{n+k})^{-1}$.

We will show that $H_{\g}$ is the inverse of $R^g$. This implies that the continuous map $R^g$ is a bijection and so is
a homeomorphism by compactness.

Let $\ss^* \in K^*_{G^*}$. By (\ref{sft07}) and (\ref{sft06})
\begin{equation}\label{sft08}
s^*_j = (S^j(\ss^*))_0 = J_{n+k}(H_{\g}(S^j(\ss^*))) = J_{n+k}(g^j(H_{\g}(\ss^*))) = R^g(H_{\g}(\ss^*))_j.
\end{equation}
That is, $\ss^* = R^g(H_{\g}(\ss^*))$, i.e. $R^g \circ H_{\g} = 1_{K^*_{G^*}}$.

Now let $x \in X$. Write $x = s_0 \bar s_0 \bar s_1 \dots $. That is, $s_0 = J_n(x)$ and
$\bar s_j = x_{n + jk} \dots x_{n + (j+1)k - 1}$.  Define $s_j = J_n(g^j(x))$. By induction and the definition of $g$
it follows that $g^j(x) = s_j \bar s_j \bar s_{j+1} \dots$ and $s_{j+1} = \g(s_j \bar s_j)$. Hence,
$J_{n+k}(g^j(x)) = s_j \bar s_j$. Thus, $\bar s_j = \bar J(R^g(x)_j)$. Finally, $s_0 = J(R^g(x)_0)$.
It follows that $H_{\g}(R^g(x)) = x$, i.e. $H_{\g} \circ R^g = 1_X$.

(c) Equation (\ref{sft03}) is the definition of the statement that both $f$ and $f_1$ are associated with $\g$.

(d) Equation (\ref{sft04}) just says that $R^g \circ Q^f = R^f$.

The metric estimates in (c) and (d) follow from (\ref{sft00}).

$\Box$ \vspace{.5cm}

Thus, the shift-like system $(X,g)$ is conjugate to the finite type subshift $(K^*_{G^*},S)$. For each $x \in X$,
$Q^f(x) \in X$ is a point whose $g$ orbit $\frac{1}{n+k+1}$ shadows the $f$ orbit of $x$.

Let $N$ be the cardinality of the alphabet $A$ (recall that $X = A^{\Z_+}$). On $X$ we will use as the background measure
$\l_0$ the $\frac{1}{N}$ Bernoulli measure.  That is, $\l_0(\langle w \rangle) = \frac{1}{N^{n(w)}}$ where $n(w)$ is the
length of the word.  So if the word $w = s_0 \bar s_0 \dots \bar s_{p-1}$ with $s_0 \in K$ and
$\bar s_0, \dots, \bar s_{p-1} \in \bar K$ then $\l_0(\langle w \rangle) = \frac{1}{N^{n + pk}}$.
For $s \in K$ we let $\l_s$ be the induced
probability measure with support $\langle s \rangle$. So
\begin{equation}\label{sft09}
\l_s ( \langle s \bar s_0 \dots \bar s_{p-1} \rangle ) \ = \   \frac{1}{N^{pk}}.
\end{equation}
Clearly,
\begin{equation}\label{sft10}
\l_0 \ = \ \Sigma_{s \in K} \l_0(s) \l_s \ = \ \frac{1}{N^{n}} \Sigma_{s \in K}  \l_s.
\end{equation}

For the special two alphabet model associated with $\g : K^* \to K$ we define the distribution data\index{distribution data}
$\nu : K^* \to [0,1]$ by
\begin{equation} \label{sft11}
 \nu(s^*)  \ = \   \frac{1}{N^{k}} \ = \ \l_s(\langle \bar s \rangle) \quad \text{with} \ s^* = s \bar s.
 \end{equation}

 By (\ref{alp13}) the measure $\mu_s$ on $K^*_{G^*}$, is given by
 \begin{equation}\label{sft12}
 \mu_s(\langle s^*_0 \dots s^*_{p-1} \rangle ) \ = \ \begin{cases} \ \frac{1}{N^{pk}} \ \text{if} \ J(s^*_0) = s, \\
 \ 0 \qquad \text{otherwise},\end{cases}
 \end{equation}
 for $s^*_0 \dots s^*_{p-1}$ a $K^*_{G^*}$ sequence.

 \begin{lem}\label{sftlem04} $(H_{\g})_*(\mu_s) = \l_s$. \end{lem}

 {\bfseries Proof:} From (\ref{sft06}) it is clear that for the $K^*_{G^*}$ sequence $s^*_0 \dots s^*_{p-1}$ with $J(s^*_0) = s$
 the homeomorphism $H_{\g}$ maps
 $\langle s^*_0 \dots s^*_{p-1} \rangle$ clopen in $K^*_{G^*}$ to  $\langle s \bar s_0 \dots \bar s_{p-1} \rangle$
 with $\bar s_j = \bar J(s^*_j)$ for $0 \leq j < p$.  By (\ref{sft09}) and (\ref{sft12}) the measures agree.

$\Box$ \vspace{.5cm}

As $(B^*,B)$ varies over the pairs of associated basic sets for $(K^*, G^*)$ and $(K, G)$, the subsets $B^*_{G^*}$ vary over the
$S$ basic sets in $K^*_{G^*}$. Let $B^*_g = H_{\g}(B^*_{G^*}) \subset X$.
Since $H_{\g}$ is a conjugacy, these are the $g$ basic sets. Furthermore, $B^*_g$ is terminal if and only if $B^*$ and $B$ are.

If $B$ is terminal then by (\ref{alp12a}) and (\ref{sft11}) the stationary distribution $v_B$ on $B$
satisfies for $s_2 \in K$:
\begin{equation}\label{sft13}
%\Sigma_{s^* \in B^* \cap \g^{-1}(s_2)} \ v_B(J(s^*)) \cdot \nu(s^*) \ = \
\frac{1}{N^k} \cdot \Sigma_{s^* \in B^* \cap \g^{-1}(s_2)} \ v_B(J(s^*))  \ = \  v_B(s_2).
\end{equation}

 \begin{theo}\label{sdtheo12a}  Let $A$ be a finite alphabet of cardinality $N$.
 Let $g$ on $X = A^{\Z_+}$ be the shift-like map associated with
  $\g : K^* \to K$ for $K = A^{[n]}$ and $K^* = A^{[n +k]}$. Let $J = J_n : \ K^* \to \ K$ be the initial word map and
 $G^* = J^{-1} \circ \g $ be the associated
relation on the finite set $K^*$.  Let $v_0$ be the positive distribution on $K$ with
$v_0(s) = \l_0(s) = \frac{1}{N^{n}}$ for $s \in K$. Let
$\l_0 = \Sigma_{s \in \ K} \ v_0(s) \l_s$, which is the $\frac{1}{N}$ Bernoulli measure, be the background measure on $X$.

Let $B^*$ and $B$ be associated terminal basic sets for $( K,G)$ and $( K^*,G^*)$.
The measure $\l_B = \Sigma_{s \in B} v_B(s) \l_s$ is an ergodic measure for $g$ with support $B^*_g = H_{\g}(B^*_G)$.

With the measures $\{ \l_B \}$ indexed by the terminal basic sets of $(K,G)$, the system $(X,g,\l_0)$ is
tractable.
\end{theo}

{\bfseries Proof:} Since $H_{\g}$ is a conjugacy which preserves the various measures, the result follows from Corollary
\ref{alpcor08}.

$\Box$

\vspace{1cm}

\bibliographystyle{amsplain}

\printindex

\end{document}